\documentclass[12pt]{article}
\usepackage[margin=1in]{geometry}
\usepackage{url}

\usepackage[T1]{fontenc}
\usepackage{fouriernc}

\usepackage{amsmath}
\usepackage{amssymb}
\usepackage{amsthm}
\usepackage{mathrsfs}

\usepackage{hyperref}
\usepackage{xcolor} 
\usepackage{graphicx}

\theoremstyle{plain}
\newtheorem{theorem}{Theorem}
\newtheorem{corollary}[theorem]{Corollary}
\newtheorem{lemma}[theorem]{Lemma}

\theoremstyle{definition}
\newtheorem{definition}[theorem]{Definition}
\newtheorem{notation}[theorem]{Notation}
\newtheorem{example}[theorem]{Example}

\newtheorem{problem}[theorem]{Problem}

\theoremstyle{remark}
\newtheorem{remark}[theorem]{Remark}

\newcommand{\pound}{\operatornamewithlimits{\circ}}
\newcommand\T{\rule{0pt}{2.6ex}}       
\newcommand\B{\rule[-1.2ex]{0pt}{0pt}} 

\title{Closed Ziv-Lempel factorization of the $m$-bonacci words}
\author{
	Marieh Jahannia\\
	School of Mathematics, Statistics and Computer Science, College of Science\\
	University of Tehran, Tehran, Iran\\
	\url{mjahannia@ut.ac.ir}\\
	\and
	Morteza Mohammad-noori\\
	School of Mathematics, Statistics and Computer Science, College of Science\\
	University of Tehran, Tehran, Iran\\
	\and
	Narad Rampersad\\
	Department of Mathematics and Statistics\\
	University of Winnipeg, Winnipeg, Canada\\
	\url{n.rampersad@uwinnipeg.ca}\\
	\and
	Manon Stipulanti\\
	Department of Mathematics\\
	University of Li\`ege, Li\`ege, Belgium\\
	\url{m.stipulanti@uliege.be}}
\date{\today}

\begin{document}
\maketitle
	
	\begin{abstract}
		A word $w$ is said to be closed if it has a proper factor $x$ which occurs exactly twice in $w$, as a prefix and as a suffix of $w$. Based on the concept of Ziv-Lempel factorization, we define the closed $z$-factorization of finite and infinite words. Then we find the closed $z$-factorization of the infinite $m$-bonacci words for all $m \geq 2$. We also classify closed prefixes of the infinite $m$-bonacci words.
	\end{abstract}
	
\bigskip
\hrule
\bigskip



\bigskip
\hrule
\bigskip
	\section{\bf Introduction}
	\vskip 0.4 true cm
	Factorization of words is an important topic in combinatorics on words, which roughly consists in breaking a given word into concatenation of other words, called factors.
	Some specific factorizations require that those factors satisfy some special properties.
	Some various types of factorizations studied in the literature are the Ziv-Lempel factorization, the Crochemore factorization, the Lyndon factorization and the grammar-based factorization \cite{Crochemore, Rytter,RytterGb,Siromony}.
	The Ziv-Lempel factorization, or $z$-factorization for short, was introduced by Ziv and Lempel for finite words \cite{Ziv}  and then was extended to infinite words \cite{BersSav}.
	This factorization has several applications in data compression \cite{Ziv1977} and text processing \cite{Kolpakov}.
	Ghareghani et al. \cite{Ghareghani} determined $z$-factorizations for standard episturmian words.
	We introduced the palindromic $z$-factorizations by requiring each factor to be a palindrome and computed this factorization for the $m$-bonacci words  \cite{jahannia}.
	In this work, based on the notion of closed words, which appeared in \cite{Bucci}, we introduce the closed $z$-factorization and apply it to the infinite Fibonacci word and then to all $m$-bonacci words, for $m>2$.
	We also characterize closed prefixes of the $m$-bonacci word $h_{\omega}$ and obtain the binary word $x=oc(h_{\omega})$ associated with  closed prefixes of $h_{\omega}$ (defined by $x_n=1$ if the prefix of length $n$ of $h_{\omega}$ is closed; otherwise, $x_n=0$). The connection of this word with the sequence of $m$-bonacci numbers then appears as a consequence.
	
	The paper is organized as follows. In Section 2, we present some notation and definitions needed in the rest of the paper. In Section 3, we study the closed $z$-factorization of the Fibonacci word. In Section 4 we study the closed $z$-factorization of the $m$-bonacci words and consider numerous properties of this factorization. Section 5 is devoted to the link between the closed and palindromic $z$-factorizations of the $m$-bonacci words.
	In Section 6, we characterize closed prefixes of the $m$-bonacci words and we give the oc-sequence of the $m$-bonacci words. Finally in Section 7 we mention some open problems.

	\section{\bf  Preliminaries}
	\vskip 0.4 true cm
	Let $A$ be a finite alphabet.
	The elements of $A^*$ are called (finite) words over $A$. We denote the empty word by $\varepsilon$ and we let $A^{+}=A^* \setminus \{ \varepsilon \}$. For every finite word $w$, we let $|w|$ denote its \emph{length}. A word $z$ is a \textit{factor} of $w\in A^*$, and we write $z \prec w$, if $w=uzv$ for some $u , v \in A^*$.
	The factor $z$ is said to be \textit{proper} if $u , v \neq \varepsilon$.
	We say that $z$ is a \textit{prefix} (resp., \textit{suffix}) of $w$, and we denote this by $z\vartriangleleft w$ (resp., $z \vartriangleright w$ ), if $u= \varepsilon$ (resp., $v= \varepsilon$).
	The set of factors of a word $w$ is denoted by $\boldsymbol{F}(w)$.
	For a factor $z$ of a word $w$, we let $|w|_z$ denote the number of occurrences of $z$ in $w$.
	We say that $z$ is a \textit{right special factor} of $w$ if $za$ and $zb$ are factors of $w$ for some distinct letters $a, b \in A$.
	
	Let $w= w_1\, w_2 \, \cdots \, w_n$ with $w_i \in A$ for all $i \in \{1,\ldots,n\}$. We let ${w}^R$ denote the \textit{reverse} of $w$, that is, ${w}^R= w_n \, \cdots \, w_2 \, w_1$. If $w={w}^R$, $w$ is called a \textit{palindrome} or a \textit{palindromic word}.
	We let $(w)^+$ denote the \textit{palindromic closure} of $w$, that is,  the shortest palindrome having $w$ as a prefix. For example, $(race)^{+}=racecar$. If $v$ is the longest palindromic suffix of $w$, say $w=uv$, then $(w)^+=uv{u}^R$. For each word $u \in A^*$, we use the notation $u^{-1}$ as below. If $w=uv$, then we let $u^{-1}w=v$ and $wv^{-1}=u$. This simply yields $(uv)^{-1}=v^{-1} u^{-1}$, consequently, $w^{-1}=w_n^{-1} \cdots w_2^{-1} w_1^{-1}$.
	
	Let ${\boldsymbol{t}}$ be an infinite word and let $w$ be a factor of ${\boldsymbol{t}}$.
	The word $v$ is said to be a \textit{return word} of $w$ if $v$ begins with an occurrence of $w$ and ends exactly just before the next occurrence of $w$ in $\boldsymbol{t}$. If $v$ is a return word of $w$, then $vw$ is said to be a \textit{complete return word} of $w$.
	The notion of return words was introduced to study primitive substitutive sequences \cite{Durand}. In \cite{Justin},
	Justin and Vuillon presented a new characterization of Sturmian words using return words. They also characterized the return words of factors of standard episturmian words.
	
	A non-empty word $x$ is called a \textit{border} of $w$ if $x$ is both a prefix and a suffix of $w$. A word $w$ is said to be \textit{closed} if it is a single letter or has a border $x$ such that it does not have any other occurrence in $w$, in other words, $|w|_{x}=2$. In this case, we call $w$ the \textit{frontier} of $x$. As an example, the word $w=mamma$ is closed, because $ma$ appears only as a prefix and a suffix of $w$. 
	The notion of closed words appeared in the study of trapezoidal words \cite{Bucci}.
	If $w$ is not closed then it is said to be {\emph{open}.

		In \cite{FiciLiptak}, Fici et al. studied words with the smallest number of closed factors. In \cite{BadkFiLi}, Badkobeh et al. showed that a length-$n$ word contains at least $n+1$ distinct closed factors and characterized those words having exactly $n+1$ closed factors. Badkobeh et al. described an efficient solution to the shortest and longest closed factorizations \cite{Badkobeh}. The shortest (resp., longest) closed factorization of a string is obtained by factorizing it into shortest (resp., longest) closed factors.
		In \cite{Luca}, A. De Luca et al. studied closed prefixes of Sturmian words and introduced the oc-sequence of a word $w$, as $oc(w)$, which is the binary sequence whose $n$-th term is $1$ if the length-$n$ prefix of $w$ is closed, or $0$ if it is open. They showed that this sequence is deeply related to the combinatorial and periodic structure of a word.
		In \cite{Bucci}, Bucci et al. studied  closed prefixes of Fibonacci words and investigated the oc-sequence of the Fibonacci word, $oc(F)$.
		Note that, for a given infinite word $\boldsymbol{t}$ and a nonempty factor $w$ of $\boldsymbol{t}$, every complete return word of $w$ is closed because it contains the factor $w$ exactly twice, once as a prefix and once as a suffix.
		
		Let $A$ be a finite alphabet. A mapping $\psi \, \colon \, A^* \to B^*$ is called a \textit{morphism} if $\psi (uv) = \psi(u) \, \psi(v)$ for all $u , v \in A^*$.
		A morphism $\psi$ is said to be \textit{prolongable} if there exists a letter $a \in A$ and a word $x \in A^*$ such that $\psi(a)=ax$ and $\psi^{i}(x) \neq \varepsilon$ for all $i \geq 0$.
		In this case, the word $\psi^{n}(a)$ is a proper prefix of $\psi^{n+1}(a)$ for all $n\ge 0$. Therefore, the infinite word $\psi^{\omega}(a)=\lim_{n\to \infty}\psi^{n}(a) $ is a fixed point of $\psi$. For every morphism $\psi \, \colon \, A^* \to B^*$ and each word $u \in A^*$, we define $\psi(u^{-1})=(\psi(u))^{-1}$. This is justified by applying $\psi$ on $uu^{-1}=u^{-1}u=\varepsilon$.
		
		A \textit{factorization} of a word consists in decomposing it into consecutive factors, which satisfy some special properties.
		Given an infinite word $\boldsymbol{w}$, the \textit{Ziv-Lempel factorization} or the \textit{z-factorization} of $\boldsymbol{w}$ is $z(\boldsymbol{w})=(z_1\,  , z_2 , \,  \ldots)$ where $z_i$ is the shortest prefix of $z_i\,z_{i+1} \cdots$ that occurs exactly once in $z_1\,z_2\,\cdots\, z_i$.
		We introduce the closed $z$-factorization $cz(\boldsymbol{w})=(z_1\,  , z_2 , \,  \ldots)$ of $\boldsymbol{w}$ by requiring that each factor $z_i$ is closed.

		\section{\bf Closed $z$-factorization of the Fibonacci word}
		
		Recall that the sequence of \textit{finite Fibonacci words} is given by $f_{-1}=1$, $f_0=0$ and $f_n = f_{n-1}\;f_{n-2}$ for all $n\geq 1$.
		Furthermore, for all $n\geq 0$, $f_n$ is the $n$-th iteration of the morphism $\sigma \colon A^* \to A^*$ on the letter 0 defined by $ \sigma (0)=01$, $ \sigma(1)=0$, that is, $f_n=\sigma^n(0)$.
		The first few Fibonacci words are given in Table~\ref{table:first-few-words-f-n}.
		The infinite Fibonacci word $f_{\omega}$ is given by $f_{\omega}=\lim\limits_{n \to \infty } f_n$.
		Equivalently we have, $f_{\omega}=\sigma^{\omega}(0)$.
		
		\begin{table}[h]
			\centering	
			\begin{tabular}{c | c c c c c c c}
				$n$&-1 & 0 & 1 & 2 & 3 & 4 & 5 \\
				\cline{1-8}
				$f_{n}$& 1 & 0 & 01 & 010 & 01001 & 01001010 &  0100101001001\\
			\end{tabular}
			\caption{The first few Fibonacci words $(f_n)_{n\ge -1}$. }
			\label{table:first-few-words-f-n}
		\end{table}
		
		The sequence of Fibonacci numbers is given by the recurrence relation $F_{n}=F_{n-1}+F_{n-2}$ for all $n \geq 1$ where $F_{-1}=1$, $F_0=1$.
		The sequence of Fibonacci words is related to the latter sequence of numbers since $F_n=|f_n|$ for all $n\ge -1$.
		
		We note that the infinite Fibonacci word belongs to the class of \emph{Sturmian words}, that is, infinite aperiodic binary words with minimal factor complexity.
		These words were presented in \cite{Morse} and are widely studied in the literature because they have several equivalent definitions and many various optimal properties, see for instance \cite[Chapter 2]{Lothare}.

		In \cite{WenWen}, Wen and Wen defined the $n$-th \textit{ singular word} $w_n$ of $f_{\omega}$ by $w_{-2}= \varepsilon$ , $w_{-1}=0$, $w_{0}=1$ and for $n \geq 1$, $w_n=a\, f_n \,b^{-1}$, where $ab \in \{01,10\}$ is the length-$2$ suffix of $f_n$.
		It is easy to see that $|w_n|=F_n$ for all $n\ge -1$.
		The first few singular words of the Fibonacci word $f_{\omega}$ are displayed in Table~\ref{table:first-few-words-w-n}.
		
		\begin{table}[h]
			\centering	
			\begin{tabular}{c |c c c c c c c c}
				$n$&-2&-1& 0 & 1 & 2 & 3 & 4 & 5 \\
				\cline{1-9}
				$w_{n}$&$\varepsilon$&0 & 1 & 00 & 101 & 00100 & 10100101 &  0010010100100\\
			\end{tabular}
			\caption{The first few singular words $(w_n)_{n \ge -2}$ of the Fibonacci word $f_\omega$. }
			\label{table:first-few-words-w-n}
		\end{table}
		
		The following lemma summarizes some properties of the singular words of the Fibonacci word $f_{\omega}$ that are useful in the following.
		
		\begin{lemma}
			\label{singular_properties_Fibo}
			$\mathrm{ \cite[ Property\, 2]{WenWen}}$
			The singular words $(w_n)_n$ of the Fibonacci word have the following properties.
			\begin{enumerate}
				\item For all  $n \geq -1$, $w_n \nprec w_{n+1}$.
				\item For all $n \geq 1$, $w_n\,=\,w_{n-2}\, w_{n-3}\, w_{n-2}$.
				\item For all $n \geq -2$, $w_n $ is a palindrome.
				\item For all $n \geq -1$, $w_n \, \nprec \, \prod_{j=-1}^{n-1} w_j$.		
			\end{enumerate}	
		\end{lemma}
		
		It is known that the infinite Fibonacci word $f_{\omega}$ can be written as the concatenation of the singular words $(w_n)_{n\ge -1}$ \cite{WenWen},
		which turns out to be the $z$-factorization of $f_{\omega}$ \cite{Fici}.
		
		\begin{lemma} $\mathrm{\cite[Proposition\, 8]{Fici}}$
			\label{prop_zfact_fibo}
			The infinite Fibonacci word $f_{\omega}$ is the concatenation of the singular words, that is,
			\begin{equation*}
				f_{\omega}=\prod_{n \geq -1} w_{n}.
			\end{equation*}
		\end{lemma}
		Our main goal in this section is to prove Theorem~\ref{ThmFibo}, which gives the closed $z$-factorization of the infinite Fibonacci word $f_{\omega}$.
		We will make use of the following two lemmas.
		
		\begin{lemma}
			\label{lem_w}
			$\mathrm{\cite[Lemma \,3]{WenWen}}$
			Let $n \geq 1$ and write $w_n w_{n+1}=u_1\, u_2\, u_3$ (or $w_{n+1} w_n= u_1\, u_2\, u_3$) with $0 < |u_1| < F_n$ and $0 < \, |u_3| \, < F_{n+1}$.
			Then $u_2$ is not a singular word.
		\end{lemma}
		
		\begin{lemma}
			\label{Lemma_w_i_closed}
			For all $n \geq -1$, $w_n$ is closed.
		\end{lemma}
		\begin{proof}	
			From Table~\ref{table:first-few-words-w-n}, the assertion can be easily verified for the values $-1 \leq  n \leq  3$. Assume that $n \geq 4$.
			It follows from Lemma~\ref{singular_properties_Fibo} that $w_n=w_{n-2}\, w_{n-3}\, w_{n-2}$. So $w_{n-2}$ is a border of $w_n$.
			It suffices to show that $w_{n-2}$ is neither a proper factor of $w_{n-2}\, w_{n-3}$ nor that of $w_{n-3}\, w_{n-2}$ . We proceed by contradiction and suppose that $w_{n-2}$ is a proper factor of $w_{n-2}\, w_{n-3}$ (the other case is similar). There exist non-empty words $u_1$ and $u_3$ over $\{0,1\}$ such that $w_{n-2}\, w_{n-3}=u_1 \, w_{n-2}\, u_3$.
			Using the fact that $|w_{n-2}|=F_{n-2}$ and $|w_{n-3}|=F_{n-3}$,
			we have $|u_1u_3|=F_{n-3}$. So $0<|u_1|<F_{n-3}<F_{n-2}$ and  $0<|u_3|<F_{n-3}$.	This contradicts Lemma~\ref{lem_w}.
		\end{proof}
		\begin{theorem}
			\label{ThmFibo}
			The closed $z$-factorization of the Fibonacci infinite word is
			\begin{equation*}
				cz(f_{\omega}) = (w_{-1} \, , w_{0} \, , \, w_{1} \,  ,\, \ldots).
			\end{equation*}
		\end{theorem}
		\begin{proof}
			{
				Based on Lemma~\ref{prop_zfact_fibo}, $z(f_{\omega}) = (w_{-1} \, , w_{0} \, , \, w_{1} \,  ,\,\ldots)$.
				By Lemma~\ref{Lemma_w_i_closed}, the factors $w_n$ are closed, which shows that this factorization is also $cz(f_{\omega})$.
			}
		\end{proof}
		\section{\bf Closed $z$-factorization of the $m$-bonacci word}
		The {\it Tribonacci word} is the most natural extension of the Fibonacci word to a three-letter alphabet and has been studied in many papers, see for instance  \cite{ArnoxRauzy,TanWen}. To describe such an extension to a finite alphabet of arbitrary size greater than 1, for every integer $m\geq 2$, we define the {\it $m$-bonacci word} as the fixed point of the morphism $\varphi_m$ given in the following definition.
		
		\begin{definition}
			\label{m-bonacci}
			For  $m\ge 2$, let $A_m= \{0 , 1 , \ldots , m-1\}$ and let $\varphi_m$ be the morphism defined by
			\begin{equation*}
				\varphi_m \colon A_m^* \to A_m^* ,\,
				0 \mapsto 01,   \ldots ,  \: (m-2) \mapsto 0 (m-1), \: 	(m-1) \mapsto 0.
			\end{equation*}
			The sequence of \emph{finite $m$-bonacci words} denoted as $(h_n^{(m)})_{n \ge 0}$, or briefly as $(h_n)_{n \ge 0}$, is given by $h_n=\varphi_m^n(0)$ for all $n\ge 0$. 	 The \emph{infinite $m$-bonacci word} $h_{\omega}$ is the fixed point of the morphism $\varphi_m$ which starts with $0$.
		\end{definition}
		
		\begin{example}
			Suppose that $m=2$. Then $A_2=\{0,1\}$ and $\varphi_2:  0 \mapsto 01,\, 1 \mapsto 0$. Also the infinite $2$-bonacci word $h_\omega^{(2)} $ is exactly the infinite Fibonacci word. Furthermore, for all	$n \geq 0$, $f_n = h_n$.
		\end{example}
		
		The first few finite $m$-bonacci words are given in Table~\ref{table:first-few-words-h-n} for some values of the parameter $m$.

		
		\begin{table}
			\centering	
			\begin{tabular}{c | c c c c c c}
				$n$& 0 & 1 & 2 & 3 & 4 & 5 \\
				\cline{1-7}
				$h_{n}^{(2)}$ & 0 & 01 & 010 & 01001 & 01001010 &   0100101001001 \\
				\cline{1-7}
				$h_{n}^{(3)}$ & 0 & 01 & 0102 & 0102010 &    0102010010201& 010201001020101020100102  \\	
				\cline{1-7}
				$h_{n}^{(4)}$ & 0 & 01 & 0102 & 01020103 & 010201030102010 &  01020103010201001020103010201 \\	
				\cline{1-7}
				$h_{n}^{(5)}$ & 0 & 01 & 0102 & 01020103 & 0102010301020104 & 0102010301020104010201030102010  \\
			\end{tabular}
			\caption{The first few words of the sequence $\big(h^{(m)}_n\big)_{n\ge 0}$  for $m\in \{2,3,4,5\}$. }
			\label{table:first-few-words-h-n}
		\end{table}
		%
		%
		%
		
		\begin{remark}
			In addition to the morphism $\varphi_m$, we define several objects related to $m$ in this section, where the parameter $m$ is deleted in the notation of most of them for the sake of clarity. These objects are  $h_n$, $u_n$,  $\mu_n$, $\psi_n$, $z_n$, $\underline{n}$, $\underline{\underline{n}}$ and $\hat{n}$. This is justifiable firstly because fixing the parameter $m$ at the beginning of the statement of each of the upcoming results removes the danger of confusion and secondly because deleting $m$ from the notation proposes a considerable simplicity and convenience in presenting formulas.
		\end{remark}
		\begin{notation}\label{undern}
			\label{modular-not}
			For every integer $n \ge 0$, let $\underline{n}:= n \mod m$.
			Note that $\underline{n} \in \{0, 1, \ldots, m-1\}$.
		\end{notation}
		
		\begin{remark}
			\label{rem_h} The words $h_n$, represented in Definition~\ref{m-bonacci}, can be defined in a recursive way as follows:
			\begin{equation}\label{defh}
				h_n =
				\begin{cases}
					0, &\text{if } n = 0;\\
					h_{n-1}\cdots h_0 n, &\text{if } 1\leq n \leq m-1;\\
					h_{n-1}\cdots h_{n-m}, &\text{if } n\geq m.
				\end{cases}
			\end{equation}
		\end{remark}
		\begin{notation}
			Let $m$ and $n$ be integers. We use the notation $m \, | \, n$ as $m$ divides $n$ and	$m \nmid n$ as $m$ does not divide $n$.
		\end{notation}	
		\begin{lemma}
			\label{lem:starthR}
	Let $n \geq 1$. The word $h_{n}^{R}$ starts with $\underline{n}0$ if $m \nmid n$; $01$ otherwise.
		\end{lemma}
		\begin{proof}
			We proceed by induction on $n\geq 1$.  The result holds true for $n=1$ as $h_{1}^{R}=10$. Assume that $n\geq 2$. There are two cases to consider according to the value of $n$.
			
	    	\textbf{Case 1}. Suppose that $n \leq m-1$. Using Equation~\eqref{defh}, we obtain $h_n^{R}=n h_{0}^{R}\cdots h_{n-1}^{R}$. As $h_0^{R}=0$, we get the result.
	    			
			\textbf{Case 2}. Suppose that $n \geq m$.  From Equation~\eqref{defh}, we deduce that $h_n^{R}=h_{n-m}^{R} \cdots h_{n-1}^{R}$. Using the induction hypothesis, $h_{n-m}^{R}$ starts with $\underline{n-m}0 = \underline{n}0$ if $m \nmid n$; $01$ otherwise. This ends the proof.\qedhere

		\end{proof}
		Similarly to the Fibonacci word which is a typical example of Sturmian words, $m$-bonacci words are typical examples of episturmian words over $A_m$.
		Episturmian words on finite alphabet naturally extend Sturmian words to larger alphabets.
		Their construction and properties can be found in \cite{DrJuPi, JusVui,juspiri} for instance.
		Inspired by the construction of standard episturmian words presented in \cite{juspiri}, 
		we need a restricted version of that construction and prior to this, some definitions are required as well.
		
		\begin{definition}
			\label{alph-psi-mu-def}
			For each letter $a\in A_m$, we define the morphism  $\psi_a^{(m)}\colon A_m^* \to A_m^*$,  $\psi_a$ for short, by $\psi_a(a)=a$, and $\psi_a(b)=ab$ for all $b \in A_m \setminus \{a\}$.
			Furthermore, we define the sequence of morphisms  $\mu_n^{(m)}\colon A_m^* \to A_m^*$,  $\mu_n$ for short, by  $\mu_0= \textrm{id}$ and $\mu_n=\psi_{\underline{0}} \circ \psi_{\underline{1}} \cdots \circ \psi_{\underline{n-1}}$ for $n>0$, where $id: A_m^{*} \rightarrow A_m^{*}$ is the identity morphism.
		\end{definition}
		
		In Table~\ref{table:first-few-mu-n}, we compute the images of letters in $A_3$ under the morphisms $\mu_0$, $\mu_1$, \ldots, $\mu_4$.
		Observe that $h_0=\mu_0(0)$, $h_1=\mu_1(1)$, $h_2=\mu_2(2)$, $h_3=\mu_3(0)$ and $h_4=\mu_4(1)$.

		\begin{table}[h]
			\centering	
			\begin{tabular}{l | ccc}
				$a \in A_3$ & 0 & 1 & 2 \\
				\hline
				$\mu_0^{(3)}(a) = id(a)$ & 0 & 1 & 2 \\
				\hline
				$\mu_1^{(3)}(a)=\psi_0(a)$ & 0 & 01 & 02\\
				\hline
				$\mu_2^{(3)}(a)=\psi_0\circ\psi_1(a)$ & 010 & 01 & 0102 \\
				\hline
				$\mu_3^{(3)}(a)=\psi_0\circ\psi_1\circ\psi_2(a)$ & 0102010 & 010201 & 0102 \\
				\hline
				$\mu_4^{(3)}(a)=\psi_0\circ\psi_1\circ\psi_2\circ\psi_0(a)$ & 0102010 & 0102010010201 & 01020100102 
			\end{tabular}
			\caption{The images of letters in $A_3$ under morphisms $\mu_i^{(3)}$, $i \in \{0,\ldots,4\}$.}
			\label{table:first-few-mu-n}
		\end{table}
		
		\begin{definition}
			\label{un-def}
			We let $\left(u_n^{(m)}\right)_{n\ge 1}$, or briefly $\left(u_n\right)_{n\ge 1}$, denote the sequence of palindromic prefixes of $h_{\omega}$ starting with $u_1=\varepsilon$ and sorted by increasing length.
		\end{definition}
		Notice that $u_1^{(m)}=\varepsilon$ and $u_2^{(m)} = 0$ for all $m\ge 2$. It is clear that $h_{\omega}=\lim_{n\rightarrow \infty} u_n$.
		In Table~\ref{table:first-few-words-u-n}, we show the first few elements of the sequence $\left(u^{(m)}_n\right)_{n\ge 1}$ for $m \in \{ 2,3,4,5\}$.
		
		\begin{table}[h]
			\centering	
			\begin{tabular}{c |  c c c c c c}
				$n$&  1 & 2 & 3 & 4 & 5&6 \\
				\cline{1-7}
				$u_{n}^{(2)}$ &$\varepsilon$ & 0 & 010 & 010010 &01001010010  & 0100101001001010010    \\
				\cline{1-7}
				$u_{n}^{(3)}$ &$\varepsilon$ & 0 & 010 & 0102010 & 01020100102010 &    010201001020101020100102010  \\	
				\cline{1-7}	
				$u_{n}^{(4)}$  &$\varepsilon$& 0 & 010 & 0102010 & 010201030102010 & 010201030102010010201030102010   \\	
				\cline{1-7}
				$u_{n}^{(5)}$  &$\varepsilon$& 0 & 010 & 0102010 &  010201030102010 & 0102010301020104010201030102010  
			\end{tabular}
			\caption{The first few words of the sequence $\left(u^{(m)}_n\right)_{n\ge 1}$ for $m\in \{2,3,4,5\}$.}
			\label{table:first-few-words-u-n}
		\end{table}
		\begin{lemma}\label{u-h-relations}
			$\mathrm{ \cite[Section\, 2.1]{juspiri}}$ The following identities hold.
			\begin{enumerate}
				\item \label{part1_uh_lem} For  all $n \geq 0$, $h_n=\mu_n(\underline{n})$.
				\item For all $n\ge 1$, $u_{n+1}=(u_n \,\, \underline{n-1})^{(+)}$.
				
				\item For all $n\ge 1$, we have $u_{n+1}=h_{n-1} u_n$.
				
				\item \label{part2_uh_lem}For all $n \geq 2$, we have $u_{n}=h_0^R\, {h_1^R} \cdots h_{n-2}^R$.
				
				\item For all $n\geq 1$,
				\begin{equation*}
					h_n =
					\begin{cases}
						u_{n+1} \, \, n, &\text{if } 1\leq n \le m-1;\\
						u_{n} u_{n-m}^{-1}, &\text{if } n \geq m.
					\end{cases}
				\end{equation*}
			\end{enumerate}
		\end{lemma}
		
		\begin{definition}
		Let $(x_n)_{n\ge 0}$ be a sequence of words and let $(\nu_n)_{(n \ge 0)}$ be a sequence of morphisms, where $x_n \in A^{*}_{m}$ and 	$\nu_n:A_m^{*} \rightarrow A_m^{*}$ for $n \ge 0$. If $a$ and $b$ be integers where $a\leq b$, then we let
	$
			{\displaystyle \pound_{i=a}^{b}} \nu_i=\nu_a \circ \nu_{a+1} \circ \cdots \circ \nu_{b}.	
	$
	Moreover, if $a>b$, then we let ${\displaystyle \pound_{i=a}^{b}} \nu_i=id$ and $\prod_{i=a}^{b}x_i=\varepsilon$.
		\end{definition}
	In the two following lemmas, we start counting the positions of letters at 1.
		\begin{lemma}
			\label{lem-pos}	
			Let $m\geq 2$ and let $k$ and $k'$ be distinct letters of $A_m$. Suppose that the word $x=x_1 \cdots x_n$ with $x_1, \ldots,x_{n-1}\in A_m\backslash \{k\}$ contains exactly $p$ occurrences of $k'$ that occur at positions $a_1,\ldots,a_p$. Moreover, let $y=\psi_k(x)$.
			\begin{enumerate}
				\item If $x_n=k$, then $|y|=2n-1$; otherwise $|y|=2n$.
				\item The word $y$ contains $n$ occurrences of the letter $k$ that appear in all odd positions and $p$ occurrences of the letter $k'$ that appear at positions $2a_1,\ldots,2a_p$.
			\end{enumerate}	
		\end{lemma}
		\begin{proof}
			{
				Let $x=x_1\cdots x_{n-1} x_n$, where $x_1,\ldots, x_{n-1} \in A_m\setminus \{k\}$ and $x_n\in A_m$. Since $y=\psi_k(x)$, we obtain
				\begin{equation}\label{defFrec}
					y =k x_1 k x_2 \cdots k x_{n-1} k y_{n},
				\end{equation}	
				where $y_n=\varepsilon$ if $x_n=k$; $y_n=x_n$ otherwise.
				Therefore, $|y|$ equals either $2n-1$ (when $x_n=k$) or $2n$ (when $x_n \neq k$). Furthermore, by Equation~\eqref{defFrec}, for each $1\leq j \leq \lfloor \frac{y}{2}\rfloor$ the equation $y_{2j}=x_j$ holds. This proves the second part of the lemma.
			}
		\end{proof}
		
		\begin{lemma}
			Let $n \ge m $. We have
			\begin{equation}
				\label{Formula:psi}
				\big({\displaystyle \pound_{i=n-m+1}^{n-1}}\psi_{\underline{i}}\big)(\underline{n-m+1}) = \prod_{j=2}^{m}\big( \, {\displaystyle \pound_{i=n-m+1}^{n-j}} \psi_{\underline{i}}\, \big)(\underline{n-j+1}).
			\end{equation}
		\end{lemma}
		\begin{proof}
			First, observe that the length of the composition in the left-hand side of  Equation~\eqref{Formula:psi} is $m-1$. B‌y applying Lemma~\ref{lem-pos} $(m-1)$ times, both sides of  Equation~\eqref{Formula:psi} equal the word $y$ of length $2^{m-1}-1$ that consists of letters $\underline{n-1}$, $\underline{n-2}$,$\ldots$, $\underline{n-m+1}$ and for each  $j \in \{1, \ldots,  m-1\}$, the letter $\underline{n-j}$ appears exactly in positions $(2t+1) 2^{m-j-1}$ with $t \in \{0, 1, \ldots , 2^{j-1}-1\}$. The result now follows.
		\end{proof}

		\begin{lemma}\label{Pn} The following identities hold.
			\begin{enumerate}
				\item For  all $n\ge m$, $\mu_{n}(\underline{n-m+1})=h_{n-1}\,h_{n-2}\, \cdots h_{n-m+1}$.
				\item For all  $1\le  n \leq  m-1$, $\mu_{n-1}(n)=u_{n}\,n$.	
				\item For  all $n \geq 1$, ${h}_{n}^R\,= \,0^{-1}\,  \varphi_m({h}_{n-1}^R) \, 0.$
				\item For   all $n \geq 2 $,  $u_{n}=\varphi_m(u_{n-1})\, 0.$
			\end{enumerate}
		\end{lemma}
		\begin{proof}{
				\begin{enumerate}
					\item The result is obtained by applying ${\displaystyle \pound_{i=0}^{n-m}} \psi_{\underline{i}}$ on both sides of Equation~\eqref{Formula:psi} and using the identities $h_j=\mu_{j}(\underline{j})=\psi_{\underline{0}} \circ \psi_{\underline{1}} \circ \cdots \circ \psi_{\underline{j-1}}(\underline{j})$  which hold for each $j \in \{n-1, \ldots , n-m+1\}$  by Part (1) of Lemma~\ref{u-h-relations}.
					\item By Part (1) of Lemma~\ref{u-h-relations}, considering $1\le  n \leq  m-1$, we obtain $h_{n-1}= \mu_{n-1}(n-1)$ and $h_n= \mu_n(n)$. Now replacing $\mu_n$ by $\mu_{n-1}  \circ \psi_{n-1}$ in the last equation easily yields
					\begin{equation}\label{aux-p}
						h_n=h_{n-1}\, \, \mu_{n-1}(n).
					\end{equation}
					On the other hand, by Parts (3) and (5) of Lemma~\ref{u-h-relations}, we get $h_n=h_{n-1} u_n \, n$. Comparing this with Equality~\eqref{aux-p} yields the result.
					\item We proceed by induction on $n \geq 1$. The result holds true for the case $n=1$ because
					$h_1^R=10=0^{-1}\varphi_m(h_0^R)0$ as $h_0^R=0$. We divide the proof into three cases according to
					the value of $n$.
					
						\textbf{Case 1}. Suppose that $n\leq  m-1$. By~\eqref{defh}, we have
						$h_n^R=nh_{0}^R\cdots h_{n-1}^R$. Using the induction hypothesis, we get 
						$h_n^R=n00^{-1}\varphi_m({h_0^R})0 \cdots 0^{-1}\varphi_m({h_{n-1}^R})0=0^{-1}\varphi_m((n-1)h_0^R\cdots h_{n-2}^R)0=0^{-1}\varphi_m(h_{n-1}^R)0$. 
						
						\textbf{Case 2}. Suppose that $n= m-1$. Using~\eqref{defh}, we have
						$h_n^R=h_{0}^R\cdots h_{n-1}^R$. By the induction hypothesis, we find that 
						$h_n^R=00^{-1}\varphi_m({h_0^R})0 \cdots 0^{-1}\varphi_m({h_{n-1}^R})0=0^{-1}\varphi_m((n-1)h_0^R\cdots h_{n-2}^R)0=0^{-1}\varphi_m(h_{n-1}^R)0$. 
						
						\textbf{Case 3}. Suppose that $n \ge m$. By~\eqref{defh}, we deduce that
						$h_n^R=h_{n-m}^R\cdots h_{n-1}^R$. By the induction hypothesis, we get 
						$h_n^R=0^{-1}\varphi_m({h_{n-m-1}^R})0 \cdots 0^{-1}\varphi_m({h_{n-1}^R})0=0^{-1}\varphi_m(h_{n-m-1}^R\cdots h_{n-2}^R)0=0^{-1}\varphi_m(h_{n-1}^R)0$. 

					\item For $n=2$, we have $u_{2}=0=\varphi_m(u_{1})\, 0$ since  $u_1=\varepsilon$. Now suppose that $n \geq 3$.
					Using Part (4) of Lemma~\ref{u-h-relations} and $h_0=0$, we obtain $u_n=0\, h_{1}^R \cdots h_{n-2}^R$. Replacing the words	 $h_{j}^R$ in the right-hand side from Part (3) of Lemma~\ref{Pn} yields $u_n=\varphi_m( h_{0}^R \, h_{1}^R \cdots h_{n-3}^R)\, 0$ whence the result is obtained by Part (4) of Lemma~\ref{u-h-relations}.\qedhere
				\end{enumerate}
			}
		\end{proof}
		We define a sequence of words $(z_n^{(m)})_{n\geq 0}$ in terms of $(h_n^{(m)})_{n\geq 0}$ that will be useful in the sequel to obtain the closed $z$-factorization of the $m$-bonacci words.		
		\begin{definition}
			\label{Def_z_m}
			We define the sequence $(z_n^{(m)})_{n\geq 0}$, denoted briefly $(z_n)_{n \geq 0}$, by $z_0=0$, $z_1=1$, $z_2=020$ and 
			\begin{enumerate}
				\item If $m=2$, then for all $n\geq 3$, $z_n=(\underline{n-3})^{-1}\,\, h_{n-3}^{R} \,\, h_{n-2}^{R}\,\underline{n-2}$.
				\item If $m \geq 3$, then
				\begin{equation*}
					z_n =
					\begin{cases}
						(n-3)^{-1}\,\, h_{n-3}^{R} \,\, h_{n-2}^R\,n\, h_{0}^R\,h_{1}^R\, \cdots h_{n-3}^R\,(n-2), &\text{if } 3 \leq n \leq m-1;\\[10pt]
						(\underline{n-3})^{-1}\,\, h_{n-3}^R \,\, h_{n-2}^R\,\, h_{n-m}^R\,\, h_{n-m+1}^R\,\, h_{n-m+2}^R \, \,\cdots\,  h_{n-3}^R\,\, \underline{n-2}, &\text{if } n \geq m.
					\end{cases}
				\end{equation*}
			\end{enumerate}	
		\end{definition}

		The recursive equation satisfied by the sequence $(|z_n|)_n$ is given in Part (1) of Lemma~\ref{Lem_start_z_m}. 
		We show in Lemma~\ref{Lem:w_z_Fibo} that for $m=2$, the words  $z_n$ are exactly   the singular words of the Fibonacci word. The case $m=3$ is studied in the following example.		
		\begin{example}
			The case $m=3$ corresponds to the Tribonacci word. For all $n\ge 3$, we have
			$$z_n^{(3)}=(\underline{n-3})^{-1}\,\, h_{n-3}^R \,\, h_{n-2}^R\,\, h_{n-3}^R\,\, \underline{n-2}.$$
		\end{example} 
		Table~\ref{table:first-few-words-z-n} shows the first few words of the sequence $(z_n^{{(m)}} )_{n \geq 0}$ for $2 \leq m  \leq 5$.
		\begin{table}[h]
			\centering	
			\begin{tabular}{c | c c c c c c}
				$n$& 0 & 1 & 2 & 3 & 4 & 5 \\
				\cline{1-7}
				$z_{n}^{(2)}$ & 0 & 1 & 00 & 101 & 00100 &  10100101\\
				\cline{1-7}
				$z_{n}^{(3)}$ & 0 & 1 & 020 & 1001 & 02010102 &  010010201020100 \\	
				\cline{1-7}
				$z_{n}^{(4)}$ & 0 & 1 & 020 & 10301 & 020100102 &  010301020101020103 \\	
				\cline{1-7}
				$z_{n}^{(5)}$ & 0 & 1 & 020 & 10301 & 0201040102 &  0103010201001020103 
			\end{tabular}
			\caption{The first few words of the sequence $(z_n^{(m)} )_{n\geq 0}$ for $m = 2,3,4,5$. }
			\label{table:first-few-words-z-n}
		\end{table}
		\begin{lemma}
			\label{Lem:w_z_Fibo}
			For all $n \ge 0$,  $w_{n-1}=z_{n}^{(2)}$.
		\end{lemma}
		\begin{proof}
			The cases $n \in \{  0, 1, 2\} $ are easily handled using Tables~\ref{table:first-few-words-w-n} and ~\ref{table:first-few-words-z-n}. Assume that $n \ge 3$.
			By  definition, we get $w_{n-1}=a\,f_{n-1}\,b^{-1}$, where $ab$ is the length-2 suffix of $f_{n-1}$.
			It can be verified by induction on $n$ that for all $n\ge 1$, $f_n$ ends with $01$ (resp., $10$) if $n$ is odd (resp., even). Therefore, we can write $w_{n-1}=\underline{n-2}\,f_{n-1}\, (\underline{n-1})^{-1}$. Moreover, using the recursive definition of $f_{n-1}$, we get $w_{n-1}=\underline{n-2}\,f_{n-2}\,f_{n-3}\,(\underline{n-1})^{-1}.$ Using Part (3) of Lemma~\ref{singular_properties_Fibo}, we have $w_{n-1}=(\underline{n-1})^{-1}\, f_{n-3}^R f_{n-2}^R\,\underline{n-2}$. On the other hand, Definition~\ref{Def_z_m} gives $z_n^{(2)}=(\underline{n-3})^{-1}\, (h_{n-3}^{(2)})^R \,(h_{n-2}^{(2)})^R\,\underline{n-2}$. By Definition~\ref{m-bonacci}, we have $h_{n}^{(2)}=f_n$ for all $n$, which completes the proof.		
		\end{proof}

		\begin{definition}\label{under2-hat-n}
			For every integer $n \geq 0$, let
			\begin{equation}\label{under-under}
				\underline{\underline{n}} =
				\begin{cases}
					\varepsilon, &\text{if } m|n;\\
					\underline{n}, &\text{otherwise.} \\
				\end{cases}
			\end{equation}
			Moreover, let
			$\hat{n}=(\underline{\underline{n}})^{-1} \underline{n}.$
			Consequently
			\begin{equation*}
				\hat{n} =
				\begin{cases}
					0, &\text{if } m|n;\\
					\varepsilon, &\text{otherwise. } \\
				\end{cases}
			\end{equation*}
			and $\hat{n}=\underline{n} (\underline{\underline{n}})^{-1}.$
		\end{definition}
		In Table~\ref{table:first-few-uu-n}, the first few values of  $\underline{\underline{n}}$ are displayed for $m \in \{ 2,3,4,5\}$.
		\begin{table}[h]
			\centering	
			\begin{tabular}{l | cccccccccc}
				& $\underline{\underline{0}}$ & $\underline{\underline{1}}$ & $\underline{\underline{2}}$ & $\underline{\underline{3}}$&$\underline{\underline{4}}$&$\underline{\underline{5}}$&$\underline{\underline{6}}$&$\underline{\underline{7}}$&$\underline{\underline{8}}$&$\underline{\underline{9}}$ \T \B \\
				\hline
				$m=2$ & $\varepsilon$ & 1 & $\varepsilon$&1& $\varepsilon$ &1& $\varepsilon$ &1& $\varepsilon$ &1\\
				$m=3$ & $\varepsilon$ & 1 & 2&$\varepsilon$&1&2&$\varepsilon$&1&2&$\varepsilon$\\
				$m=4$ & $\varepsilon$ & 1 & 2 & 3& $\varepsilon$&1& 2 & 3& $\varepsilon$&1 \\
				$m=5$ & $\varepsilon$ & 1 & 2 & 3 &4&$\varepsilon$& 1 & 2& 3 &4
			\end{tabular}
			\caption{The first few values of $\underline{\underline{n}}$ for $m \in \{2,3,4,5\}$.}
			\label{table:first-few-uu-n}
		\end{table}
		\begin{lemma}
			\label{varphi}
			For all $n \geq 0$,
			 $\varphi_m(\underline{n}) = 0\, \underline{\underline{n+1}}$.
		\end{lemma}
		\begin{proof}
			The proof follows immediately from the definition of the morphism $\varphi_m$ and Definition~\ref{under2-hat-n}.
		\end{proof}
		In the first part of the following lemma, we study the length of the word $z_n$ and in the second part, we find the first and the last letters of $z_n$.
		\begin{lemma}\label{Lem_start_z_m} 
			\begin{enumerate}
				\item For all $n \geq m+1$, 
				$ |z_n|\, = \, |z_{n-1}|\, + \, |z_{n-2}|\, +\,\cdots\,+\, |z_{n-m}|.
				$
				\item Let $n \geq 2$. The word $z_n$ ends with $\underline{n-2}$. Moreover, $z_n$ starts with $1$ if  $m | n-3$; $0$ otherwise.	
			\end{enumerate}
		\end{lemma}
		\begin{proof}
			\begin{enumerate}
				\item 
				The case $m=2$ follows easily by induction on $n$ and Definition~\ref{Def_z_m}. Suppose that $m\geq 3$. We proceed by induction on $n$. Using Definition~\ref{Def_z_m} and Equation~\ref{defh}, we get
				\begin{equation}
					\label{length_z}
					|z_n|=
					\begin{cases}
						|h_{n-3}|+ 2|h_{n-2}|,&\text{if }  3 \leq n \leq  m-1;\\
						|h_{n-3}| +|h_{n-2}|+|h_{n-m}|+|h_{n-m+1}|+\cdots+|h_{n-3}|,&\text{if } n \geq m.
					\end{cases}     
				\end{equation} 
				For the base case  $n=m+1$ we have,
				$|z_{m+1}|= |h_{m-2}|+|h_{m-1}|+|h_{1}|+|h_{2}|+\cdots+|h_{m-2}|.
				$
				From Equation~\eqref{defh}, we know that 
				$
				|h_{m-2}|=           |h_{0}|+|h_{1}|+\cdots+|h_{m-3}|+1.
				$
				So, we get
				$|z_{m+1}|= |h_{m-1}|+3 |h_{m-2}|-2$ since $|h_{0}^{R}|=1$.  
				By using~\eqref{length_z} several times, we obtain
				$
				|z_{m}|+|z_{m-1}|+\cdots + |z_{1}|= |h_{m-1}|+3 |h_{m-2}|-2=|z_{m+1}|
				$, as desired.\\
				Now, suppose that $n\geq m+2$. By using~\eqref{length_z} several times, we have
				$$
				|z_{n-1}|+|z_{n-2}|+\cdots + |z_{n-m}|=   |h_{n-3}| +|h_{n-2}|+|h_{n-m}|+|h_{n-m+1}|+\cdots+|h_{n-3}|=|z_{n}|.
				$$
				\item Using Definition~\ref{Def_z_m},  the result holds true for the case $n=2$ as $z_2 = 010$. Assume
				that $n \geq 3$.
				By Definition~\ref{Def_z_m} again, $z_n$ starts with $(\underline{n-3})^{-1} h_{n-3}^{R}$ and ends with $\underline{n-2}$.
				If $m|n - 3$, then Lemma~\ref{lem:starthR} says that 
				$h_{n-3}^{R}$ starts with $01 =\underline{ n - 3}\,1$ which ends	the proof.
\qedhere
			\end{enumerate}
		\end{proof}
		
		In the following remark, we rewrite the definition of the sequence $z_n$ in a slightly different form. This helps us to prove Lemma~\ref{Lem_z_m_phi} which gives a characterization of $z_n$ in terms of the morphism $\varphi_m$ and the word $z_{n-1}$.

		\begin{remark}
			\label{Rem_z_m} For all  $n\geq 3$,
			the word $z_n$ can be written as 
			$z_n =\alpha_{n}\,\, h_{n-3}^R \,\, h_{n-2}^R\,\beta_{n}\,z'_{n}\,\gamma_{n},$
			where
			$$ z'_{n}=\prod_{i={(n-m)_*}}^{(n-3)} h_i^R,\quad (j)_*:=\max\{0,j\}, \quad		
			(\alpha_{n}, \beta_{n}, \gamma_{n}) =
			\begin{cases}
				\left((n-3)^{-1}, n, n-2\right),&\text{if }  n \leq m-1;\\
				(({\underline{n-3}})^{-1}, \varepsilon, \underline{n-2}),&\text{if } n \geq m.
			\end{cases}
			$$				
			
		\end{remark}
		
		\begin{lemma}
			\label{cor:phi:h}
			For any sequence $(i_1, \ldots, i_k)$  of positive integers,  we have
			$$
			\varphi_m \big( \prod_{j=1}^{k} h_{i_j}^R \big) = 0\, \big(\prod_{j=1}^{k}h_{i_j+1}^R \big)\,0^{-1}.
			$$
			Consequently, for all $n\geq m + 1 $,
			$$
			\varphi_m \big( \prod_{i=n-m-1}^{n} h_{i}^{R} \big) = 0\, \big(\prod_{i=n-m}^{n+1}h_{i}^{R} \big)\,0^{-1}.
			$$ 	
		\end{lemma}
		\begin{proof}
			To prove the first part, from Part (3) of Lemma~\ref{Pn}, we  find that 
			$$
			\varphi_m ( h_{i_1}^R h_{i_2}^R\cdots h_{i_k}^R) 
			=0 h_{i_1+1}^R \,0^{-1}\, 0\, h_{i_2+1}^R\, 0^{-1}\, \cdots\, 0\, h_{i_k+1}^R\, 0^{-1}.
			$$ 	
			The second part is obtained from the first one.
		\end{proof}
		
		\begin{lemma}
			\label{Lem_z_m_phi}
			For all $n \geq 2$, $z_n=(\widehat{n-3})^{-1} \, \varphi_m(z_{n-1})\, \widehat{n-2}. $
		\end{lemma}
		\begin{proof}
			The proof is by induction on $n$. Assume that the result holds true up to $n-1$.	If $n =2$, then $z_2=\varphi_m(z_1)0=(\widehat{2-3})^{-1}\varphi_m(z_1)\,\widehat{2-2}.$
			If $n=3$, then
			$z_3=0^{-1}\varphi_m(z_2)=(\widehat{3-3})^{-1}\varphi_m(z_2)\,\widehat{3-2}.$
			Now suppose that $n\ge 4$.  Let
			$B_n= (\widehat{n-3})^{-1}  \varphi_m(z_{n-1}) \widehat{n-2}$. We divide the proof into two cases according to the value of $n$.
		
			\textbf{Case 1}. Suppose that $n \leq m$. Using the induction hypothesis and Remark~\ref{Rem_z_m}, we have
				\begin{align*}
					B_n	&= (\widehat{n-3})^{-1} \varphi_m\big((n-4)^{-1} h_{n-4}^R  h_{n-3}^R (n-1) \prod_{i=0}^{n-4} h_i^R\, (n-3)\big)\,\widehat{n-2}.	
				\end{align*}
			\noindent From Part (3) of Lemma~\ref{Pn} and Lemmas~\ref{varphi} and~\ref{cor:phi:h}, we get
				$$
				B_n=(\widehat{n-3})^{-1} (0 (n-3))^{-1}  0 h_{n-3}^R \, h_{n-2}^R \,a \, \prod_{i=0}^{n-3}h_{i}^R 0^{-1} 0 (n-2)\,\, \widehat{n-2},
				$$
				where $a=\varepsilon$ if $n=m$; $a=n$ otherwise.
				As a consequence, we find that
				$$
				B_n=(\underline{n-3})^{-1} h_{n-3}^R \, h_{n-2}^R\, a \, \prod_{i=0}^{n-3}h_{i}^R\underline{n-2}.
				$$
				The result can be easily deduced using Remark~\ref{Rem_z_m}.
					
				\textbf{Case 2}. Assume that $n \geq m+1$. From the induction hypothesis and then Remark~\ref{Rem_z_m}, we get
				\begin{equation*}
					B_n	= (\widehat{n-3})^{-1} \varphi_m\big((\underline{n-4})^{-1} h_{n-4}^R  h_{n-3}^R \prod_{i=n-m-1}^{n-4} h_i^R (\underline{n-3})\big)\;  \widehat{n-2}.	
				\end{equation*}
				Part (3) of Lemma~\ref{Pn} and Lemmas~\ref{varphi} and~\ref{cor:phi:h} implies that
				$$		B_n=(\widehat{n-3})^{-1} (\,0\, (\underline{\underline{n-3}}))^{-1}  \,0\, h_{n-3}^R \, h_{n-2}^R \, \prod_{i=n-m}^{n-3}h_{i}^R \,0^{-1} \,0\, (\underline{\underline{n-2}})\,\, \widehat{n-2}.
				$$
				Using Definition~\ref{under2-hat-n}, we obtain
				$$	B_n=(\underline{n-3})^{-1} h_{n-3}^R \, h_{n-2}^R   \prod_{i=n-m}^{n-3}h_{i}^R\,\underline{n-2}.
				$$
				Consequently, Remark~\ref{Rem_z_m} gives that
				$
				B_n= z_n.
				$\qedhere	
		\end{proof}
	\begin{definition}
		For   $n \geq 0$, let $P_n=\prod_{k=0}^{n-1} z_{k}$.
	\end{definition}
		\begin{lemma}
			\label{Lem_Pn}
			  We have $P_0=\varepsilon$, $P_1=0$, $P_2=01$ and for all $n \geq 3$, 
			\begin{equation*}
				P_n=\varphi_m(P_{n-1})\,\,\widehat{n-3}.		
			\end{equation*}		
		\end{lemma}
		
		\begin{proof}	
			The cases $n\in\{0,1,2,3\}$ are trivially verified. Assume that $n\ge 4$. By Lemma~\ref{Lem_z_m_phi}, we get
			\begin{align*}
				P_n&=z_0\,z_1\,\cdots \, z_{n-1}=
				0 \cdot 1\cdot (\widehat{-1})^{-1}\varphi_m(z_1)  \widehat{0}\cdot \,(\widehat{0})^{-1} \varphi_m(z_2) \widehat{1} \cdot \, \cdots \, \cdot (\widehat{n-4})^{-1}\, \varphi_m(z_{n-2})\,\,\widehat{n-3} \\
				&= \varphi_m( z_0 \,z_1\, \cdots z_{n-2})\,\, \widehat{n-3}=\varphi_m(P_{n-1})\,\, \widehat{n-3}.\qedhere
			\end{align*}
		\end{proof}
		\begin{lemma}
			\label{Lem_factorization_z_m}
			We have the following factorization for the infinite $m$-bonacci word.
			\begin{equation*}
				h_{\omega} \,=\, \prod_{n \geq 0} {z}_{n}.
			\end{equation*}
		\end{lemma}
		
		\begin{proof}
			We need to prove the following two items.
			\begin{enumerate}
				\item For all $n\ge 1$, $|P_n| > |P_{n-1}|$,
				\item $(P_n)_{n\ge 0}$ is a sequence of prefixes of $h_{\omega}$.
			\end{enumerate}

			To prove (1), note that $|P_n| = |P_{n-1}| + |z_{n-1}|$ for all $n\ge 1$. This yields $|P_{n}| > |P_{n-1}|$ since $|z_{n-1}|>0$.
			Let us prove (2). We proceed by induction on $n$.
			The $m$-bonacci word $h_{\omega}$ starts with $01$. Therefore, it is clear that $P_n$ is a prefix of $h_{\omega}$ for $n \in \{ 0,1,2\}$.
			Now suppose that $n\ge 3$ and $P_{n-1}$ is a prefix of $h_{\omega}$. Using Lemma~\ref{Lem_Pn}, we have  $P_n= \varphi_m(P_{n-1}) \widehat{n-3}$.
			The proof is divided into two cases according to whether $m$ divides $n-3$ or not.

			\textbf{Case 1}.  Suppose that $m|n-3$. Thus, $\widehat{n-3}=0$ and then $P_n= \varphi_m(P_{n-1}) 0$.	
			By the induction hypothesis, there exists a letter $a\in A_m$ and an infinite word $ \boldsymbol{z}$ over $A_m$ such that $h_\omega = P_{n-1}a \boldsymbol{z}$.
			Since $h_{\omega}$ is the fixed point of $\varphi_m$, we know that
			$
			h_{\omega} = \varphi_m(h_{\omega}) = \varphi_m(P_{n-1}\,a\boldsymbol{z}).
			$
			Since $\varphi_m(a)=0b$ with $b \in\{\varepsilon,  0,1,\ldots,m-1\}$, we get
			$
			h_{\omega}=\varphi_m(P_{n-1})\, 0\, b\, \varphi_m(\boldsymbol{z}) = P_n\, b\, \varphi_m(\boldsymbol{z}),
			$
			showing that $P_n$ is also a prefix of $h_{\omega}$.
			
			\textbf{Case 2}. Suppose that $m \nmid\,n-3$. Therefore, $\widehat{n-3}=\varepsilon$ and then, $P_n= \varphi_m(P_{n-1})$.	
			Thus $\varphi_m(P_{n-1})=P_{n}$ is a prefix of $\varphi_m(h_{\omega})=h_{\omega}$, as required. \qedhere	
		\end{proof}
		
		Justin and Vuillon \cite{Justin} studied the return words of factors of standard episturmian words and their occurrences.
		In order to prove the main theorem of this section, we need to mention some useful lemmas.  Recall from Definition~\ref{un-def} that $u_n$ is the sequence of palindromic prefixes of $h_\omega$.
		\begin{lemma}
		$\mathrm{	\cite[Corollary\, 4.1]{Justin}}$
			\label{Lem_Justin}
			Let $v \in A_m^*$ be any finite factor of $h_{\omega}$.
			Let $j(v)\ge 1$ be such that  $u_{j(v)}$ is the shortest palindromic prefix of $h_{\omega}$ which contains $v$ as a factor, say  $u_{j(v)} = fvg$ with $f, g \in A_m^{*}$. Then, $y$ is a return word of $v$ if and only if $fyf^{-1}$ is a return word of $u_{j(v)}$. Moreover, the return words of the palindromic prefix $u_{{j(v)}+1}$ are $\mu_{j(v)}(i)$ for all $i \in A_m$.
		\end{lemma}

		The following definition, which will be useful in this paper, is mentioned in \cite{Lothare}.
		\begin{definition}
			Given an alphabet $A$, a set $X \subset A^+$ of non-empty words is a code on $A$ if every word $w \in A^*$ has at most one factorization using words of $X$.	 
		\end{definition}
		\begin{lemma}
			$\mathrm{\cite[Lemma\,15]{jahannia}}$
			\label{lem:code}
			The set $\{ 0\,\underline{\underline{i}} \;|\; 1 \leq i \leq m\}=\{01,02,\ldots,0(m-1),0\}$ of non-empty words is a code on the alphabet $A_m$.
		\end{lemma}
		\begin{lemma}
		$\mathrm{	\cite[Lemma \,17]{jahannia}}$
			\label{Lem_varphi_factor}
			Let $x,y \in A_m^*$ be two finite words.
			\begin{enumerate}
				\item If $\varphi_m(x)0$ is a factor of $\varphi_m(y)0$, then $x$ is a factor of $y$.
				\item If $\varphi_m(x)$ is a factor of $\varphi_m(y)$ and $x$ does not end with the letter $m-1$, then $x$ is a factor of $y$.
			\end{enumerate}
		\end{lemma}
		
		\begin{lemma}
			\label{Lem_t_accurred_in_u_}
			For $n \geq 3$,  $(\underline{n-3}) ^{-1}\, h_{n-3}^R\, \underline{n-2}$  is not a factor of $u_{n-1}$.
		\end{lemma}
		\begin{proof}
			The proof is by induction on $n$. The statement can be readily verified for $n=3$.  Now assume that $n \geq 4$ and that the assertion holds for values less than $n$. We proceed by contradiction. Suppose that 
			$(\underline{n-3}) ^{-1}\, h_{n-3}^R\,\,\underline{n-2}\, \prec\, u_{n-1}.$
			Using Parts (3) and  (4) of Lemma~\ref{Pn}, we have
			$
			(\underline{n-3})^{-1}\,0^{-1}\,\varphi_m( h_{n-4}^R)\,0\,\,\underline{n-2} \,\prec \, \varphi_m(u_{n-2})\,0.
			$
			So we get
			\begin{equation*}
				(\underline{n-3})^{-1}\,0^{-1}\,\varphi_m\big(
				(\underline{n-4})\,(\underline{n-4})^{-1}
				h_{n-4}^R  \, (\underline{n-3})(\underline{n-3})^{-1}\big)\,0\,\,\underline{n-2}\, \prec \, \varphi_m(u_{n-2})\,0.
			\end{equation*}
			From Lemma~\ref{varphi} and then Definition~\ref{under2-hat-n}, we find that
			\begin{equation}
				\label{eq:not_occured_u}
				(\widehat{n-3})^{-1}\varphi_m\big(
				(\underline{n-4})^{-1} \,h_{n-4}^R \,(\underline{n-3})\big)\,\widehat{n-2}\, \prec \, \varphi_m(u_{n-2})\,0.
			\end{equation}
			We consider three cases according to whether $m$ divides $n-3$ or $n-2$ or none of them.
			
			\textbf{Case 1}. Suppose that $m  \nmid \,n-3$ and $m  \nmid \,n-2$. Replacing $\widehat{n-3}= \widehat{n-2}=\varepsilon$ into~\eqref{eq:not_occured_u}, we get
				$\varphi_m\big( \,(\underline{n-4})^{-1} h_{n-4}^R\, \underline{n-3}\big)\,\, \prec \, \varphi_m(u_{n-2})\,0.$
				Since $\varphi_m\big( \,(\underline{n-4})^{-1} h_{n-4}^R\, \underline{n-3}\big)$ ends with $\underline{\underline{n-2}}\neq 0$, we find that
				$\varphi_m\big( \,(\underline{n-4})^{-1} h_{n-4}^R\, \underline{n-3}\big)\,\, \prec \, \varphi_m(u_{n-2})\,.$
				Since $ \,(\underline{n-4})^{-1} h_{n-4}^R\, \underline{n-3}$ ends with $\underline{n-3} \neq m-1$, we deduce from Part (2) of Lemma~\ref{Lem_varphi_factor} that  $(\underline{n-4})^{-1}\,                                  h_{n-4}^R\,\, \underline{n-3}\, \prec\, u_{n-2}$,  contradicting the induction hypothesis.	
				
				\textbf{Case 2}. Suppose that $m  \mid \,n-3$ and $m  \nmid \,n-2$. Substituting $\widehat{n-3}=0$ and $\widehat{n-2}=\varepsilon$ into~\eqref{eq:not_occured_u}, we have 	
				$0^{-1}\,\varphi_m\big( \,(\underline{n-4})^{-1}h_{n-4}^R\,\,\underline{n-3}\big)\, \prec \, \varphi_m(u_{n-2})\,0$. By  
				Lemma~\ref{lem:starthR}, $0^{-1}\,\varphi_m\big( \,(\underline{n-4})^{-1} h_{n-4}^R\,\,\underline{n-3}\big)$ starts with $1$. From Lemma~\ref{lem:starthR} and Part (4) of Lemma~\ref{u-h-relations},  $\varphi_m(u_{n-2})$ starts with $0$. 	Hence,
				$\varphi_m\big( (\underline{n-4})^{-1} h_{n-4}^R\underline{n-3}\big)$ is a factor of $\varphi_m(u_{n-2})\,0.$
				As in Case (1), we reach a contradiction since $\underline{\underline{n-2}}\neq 0$ and $\underline{n-3}\neq m-1$.
				
				\textbf{Case 3}.  Suppose that $m  \nmid \,n-3$ and $m  \mid \,n-2$.
				By replacing $\widehat{n-3}=\varepsilon$ and $\widehat{n-2}=0$ into~\eqref{eq:not_occured_u}, we find that
				$\varphi_m\big((\underline{n-4})^{-1}	h_{n-4}^R\,\underline{n-3}\big)\,\,0\, \prec \, \varphi_m(u_{n-2})\,\,0.$
				From Part (1) of  Lemma~\ref{Lem_varphi_factor},
				$(\underline{n-4})^{-1}\,                                  h_{n-4}^R\,\, \underline{n-3}\, \prec \, u_{n-2}$, contradicting the induction hypothesis.\qedhere
		\end{proof}
		\begin{lemma}
			\label{Lem_zn_notAppear_zn+1}
			For all $n \geq 0$, $z_n$ is not a factor of $z_{n+1}$.
		\end{lemma}
		\begin{proof}If  $m=2$, then the result follows from Part (1) of Lemma~\ref{singular_properties_Fibo} and Lemma~\ref{Lem:w_z_Fibo}. We prove it for $m \geq 3$. 
			The result is obviously true  for $n \in \{ 0,1,2\}$, see Table~\ref{table:first-few-words-z-n}. Let $n\geq 3$. By induction on $n$, we assume that the result holds true up to $n-1$ and we show that it is still true for $n$. By contradiction, suppose that $z_{n} \,\prec\, \,z_{n+1}\,.$
			We obtain from Lemma~\ref{Lem_z_m_phi} that
			$
			(\widehat{n-3})^{-1} \, \varphi_m(z_{n-1}) \,\, \widehat{n-2} \, \prec \,
			(\widehat{n-2})^{-1} \,\, \varphi_m(z_{n}) \,\, \widehat{n-1}.
			$
			So, we get
			\begin{equation}
				\label{eq:lem_zn}
				(\widehat{n-3})^{-1} \, \varphi_m(z_{n-1}) \,\,  \widehat{n-2} \,\prec
				\, \varphi_m(z_{n}) \,\, \widehat{n-1}.
			\end{equation}
			We divide the proof into two cases according to whether $m$ divides $n-3$ or not.
			
			\textbf{Case 1}. Suppose that $m\nmid \,n-3$. Substituting $\widehat{n-3}=\varepsilon$ into~\eqref{eq:lem_zn}, we find that
				\begin{equation}
					\label{eq:lem_zn_case}
					\varphi_m(z_{n-1}) \,\,\widehat{n-2} \,\prec
					\, \varphi_m(z_{n}) \,\, \widehat{n-1}.
				\end{equation}
				Now we consider three cases according to  whether $m$ divides $n-1$ or $n-2$ or none of them.
				
				\textbf{Case 1-1}. Assume that $m\nmid \,n-1$ and $m\nmid \,n-2$.
					Plugging $\widehat{n-1}=\widehat{n-2}=\varepsilon$ into~\eqref{eq:lem_zn_case}, we have  $\varphi_m(z_{n-1})  \, \prec \, \varphi_m(z_{n})$. Part (2) of Lemma~\ref{Lem_start_z_m} implies that $z_{n-1}$ ends with $\underline{n-3}$. 
					Since $m\nmid \,n-2$ implies $\underline{n-2} \neq m$ and thus $\underline{n-3}\neq m-1$, using Part (2) of Lemma~\ref{Lem_varphi_factor}, $z_{n-1}$ is a factor of $z_n$, contradicting the induction hypothesis.
				
					\textbf{Case 1-2}. Suppose that $m | \,n-1$ and $m\nmid \,n-2$.
					Replacing  $\widehat{n-1}= 0$ and  $\widehat{n-2}=\varepsilon$
					into~\eqref{eq:lem_zn_case}, we find that
					$\varphi_m(z_{n-1}) \, \prec \, \varphi_m(z_{n})\,0$.
					Using Part (2) of Lemma~\ref{Lem_start_z_m},
					$z_{n-1}$ ends with $\underline{n-3}$. It follows that  $\varphi_m(z_{n-1})$ ends with $0\, \underline{\underline{n-2}}$.
					In conclusion, $\varphi_m(z_{n-1})\,  \prec \, \varphi_m(z_{n}) .$ As in Case (1-1), we reach a contradiction.
					
					\textbf{Case 1-3}. Suppose that $m \nmid \,n-1$ and $m\mid \,n-2$. Plugging $\widehat{n-1}= \varepsilon$ and  $\widehat{n-2}= 0$ into~\eqref{eq:lem_zn_case}, we get
					$\varphi_m(z_{n-1})\, 0\, \prec \, \varphi_m(z_{n})\,$. We then have,  $\varphi_m(z_{n-1})\,\, 0 \prec \, \varphi_m(z_{n})\,\,0$.  Part (1) of Lemma~\ref{Lem_varphi_factor} gives, $z_{n-1} \,\prec \,z_n$ which contradicts the induction hypothesis.
				
				\textbf{Case 2}.  Suppose that $m |\,n-3$. By replacing $\widehat{n-3}=0$ into~\eqref{eq:lem_zn}, we get  $0^{-1}\,\varphi_m(z_{n-1}) \, \widehat{n-2}\,\prec \,\varphi_m(z_n)$.
				We deduce from Part (2) of Lemma~\ref{Lem_start_z_m} that $z_n$ (resp., $z_{n-1}$) starts with 1 (resp., 0).  So $0^{-1}\varphi_m(z_{n-1})$ (resp., $\varphi_m(z_n)$) starts with 1 (resp., 0). Since each occurrence of $1$ in $\varphi_m(z_n)$ is preceded by a $0$, we conclude that 
				$
				\label{eq:lem_zn_case2}
				\varphi_m(z_{n-1})\,\, \widehat{n-2}\, \prec\, \varphi_m(z_n).
				$ and thus,
				$\varphi_m(z_{n-1})\, \prec \,\varphi_m(z_n).$
				As in Case (1-1), we reach a contradiction since $\underline{n-3} \neq m-1$.	 \qedhere
				%
		\end{proof}
		\begin{lemma}
			\label{Lem_zn_notAppearBefore}
			Let $n \geq  1$ and $\gamma_n$ be the last letter of $z_n$. The word $z_n$ is not a factor of $z_{n-1}\,z_{n}\,\gamma_{n}^{-1}$.
		\end{lemma}
		\begin{proof} To prove the result, we proceed by induction on $n$.
			If $n=1$, then $z_1=1$ is not a factor of $z_0 z_1 \gamma_1^{-1} = z_0=0$. If $n=2$, then the word
			$z_1 z_2 \gamma_2^{-1}$ does not contain $z_2$.	Now assume that $n \geq 3$ and that the result is true for all values less than $n$. We proceed by contradiction and suppose that 
			$	z_n \,\prec\, z_{n-1}\,z_{n} \gamma_{n}^{-1}.	$
			By Part (2) of Lemma~\ref{Lem_start_z_m}, $\gamma_{n}=\underline{n-2}$. Lemma~\ref{Lem_z_m_phi} implies that
			\begin{equation*}
				(\widehat{n-3})^{-1} \, \varphi_m(z_{n-1})\,\, \widehat{n-2} \,\prec\,
				(\widehat{n-4})^{-1} \, \varphi_m(z_{n-2})\,\, \widehat{n-3}\,\, (\widehat{n-3})^{-1} \, \varphi_m(z_{n-1})\,\, \widehat{n-2}\,\, (\underline{n-2})^{-1}.
			\end{equation*}
			So, we have
			$
			(\widehat{n-3})^{-1} \, \varphi_m(z_{n-1})\,\, \widehat{n-2} \, \prec
			\, \varphi_m(z_{n-2}\, z_{n-1}) \,\, (\underline{\underline{n-2}})^{-1}.
			$
			Therefore, 
			\begin{equation*}
				(\widehat{n-3})^{-1} \, \varphi_m(z_{n-1})\,\, \widehat{n-2}\, \prec
				\, \varphi_m(z_{n-2}\, z_{n-1}(\underline{n-3})^{-1} \, \underline{n-3})\,\,(\underline{\underline{n-2}})^{-1}.
			\end{equation*}
			Using Definition~\ref{under2-hat-n}, we find that
			\begin{equation}
				\label{eq:lem:notAppear}
				(\widehat{n-3})^{-1} \, \varphi_m(z_{n-1})\,\,\widehat{n-2}  \,\prec \,  \, \varphi_m(z_{n-2}\, z_{n-1} \, (\underline{n-3})^{-1})\,\,0.
			\end{equation}
			We divide the proof into three cases according to  whether $m$ divides $n-3$ or $n-2$ or none of them.

			\textbf{Case 1}. Assume that $m \nmid n-3$ and $m \nmid n-2$.
			Plugging 
			$\widehat{n-3}=\widehat{n-2}=\varepsilon$
			into~\eqref{eq:lem:notAppear}, we find that
			$\varphi_m(z_{n-1})\,  \prec \,  \, \varphi_m\big(z_{n-2}\, z_{n-1} \, (\underline{n-3})^ {-1}\big)\,\,0.$
			From Part (2) of Lemma~\ref{Lem_start_z_m}, $z_{n-1}$ ends with $\underline{n-3}$. Thus, $\varphi_m(z_{n-1})$ ends with $0\, \underline{\underline{n-2}}$. Using Definition~\ref{under2-hat-n}, $\underline{\underline{n-2}}\neq 0$. It follows that
			$\varphi_m(z_{n-1})\, \prec \,  \, \varphi_m\big(z_{n-2}\, z_{n-1} \, (\underline{n-3})^{-1}\big)$.
			 As $z_{n-1}$ ends with $\underline{n-3}$,  Part (2) of Lemma~\ref{Lem_varphi_factor} tells us that $z_{n-1}\, \prec \,z_{n-2}\,z_{n-1} \,(\underline{n-3})^{-1} $, since $m\nmid \,n-2$ implies $\underline{n-2} \neq m$ and thus $\underline{n-3}\neq m-1$.	 This contradicts the induction hypothesis.
		 	
			\textbf{Case 2}.  Assume that $m \mid n-3$ and $m \nmid n-2$.
			By replacing 
			$\widehat{n-3}= 0$ and $\widehat{n-2}=\varepsilon$
			into~\eqref{eq:lem:notAppear}, we have
			$0^{-1}\,\varphi_m(z_{n-1})\,  \prec \,  \, \varphi_m\big(z_{n-2}\, z_{n-1} \, (\underline{n-3})^ {-1}\big)\,\,0.$
			By Part (2) of Lemma~\ref{Lem_start_z_m}, $z_{n-1}$ and $z_{n-2}z_{n-1}$ start with 1.  So $0^{-1}\varphi_m(z_{n-1})$ starts with 2 and $\varphi_m(z_{n-2}z_{n-1})$ starts with 0.
			We conclude that $\varphi_m(z_{n-1})\,\prec \,\varphi_m\big(z_{n-2}z_{n-1}(\underline{n-3})^{-1}\big)\,\,0$.
			By Part (2) of Lemma~\ref{Lem_start_z_m},
			$\varphi_m(z_{n-1})$ ends with $01$ since $\underline{n-3}=0$. So,
			$\varphi_m(z_{n-1})\, \prec \,  \, \varphi_m\big(z_{n-2}\, z_{n-1} \, (\underline{n-3})^{-1}\big).$ As in Case (1), we reach a contradiction.
			
			\textbf{Case 3}.  Assume that $m \nmid n-3$ and $m \mid n-2$.
			Plugging  $\widehat{n-3}= \varepsilon$ and $\widehat{n-2}=0$ into~\eqref{eq:lem:notAppear}, we have
			$\,\varphi_m(z_{n-1})\,\,0\,  \prec \,  \, \varphi_m\big(z_{n-2}\, z_{n-1} \, (\underline{n-3})^ {-1}\big)\,\,0.$
			Part (1) of Lemma~\ref{Lem_varphi_factor} implies that  $z_{n-1}\, \prec\, z_{n-2}\, z_{n-1} \, (\underline{n-3})^{-1}$, contradicting the induction hypothesis.
			\qedhere	
		\end{proof}
		\begin{lemma}
			\label{Prop_z_Not_Appear_Before}
			For all $n \geq 1$, $z_n$ is not a factor of $P_n$.
		\end{lemma}
		\begin{proof}
			Using Part (4) of Lemma~\ref{singular_properties_Fibo} and Lemma~\ref{Lem:w_z_Fibo}, the assertion is true for the case $m=2$. Let us suppose that $m \geq 3$. The proof is by induction on $n$. It can be easily checked for $n\in\{1,2\}$ using Table~\ref{table:first-few-words-z-n}. Now suppose that $z_j$ is not a factor of $P_j$ for all $3 \leq j \leq n-1$. We show it is still true for $j=n$. Arguing by contradiction, assume that $z_{n} \,\prec \,P_{n}\,.$
			From Lemmas~\ref{Lem_z_m_phi} and~\ref{Lem_Pn}, 
			\begin{equation}
				\label{eq:lem_zn_not_appear}
				(\widehat{n-3})^{-1} \, \varphi_m(z_{n-1}) \, \widehat{n-2} \, \prec
				\, \varphi_m(P_{n-1}) \,\, \widehat{n-3}.
			\end{equation}
			The proof is divided into three cases according to whether  $m$ divides $n-3$ or $n-2$ or none of them.
			
			\textbf{Case 1}. Assume that $m \nmid n-3$ and $m \nmid n-2$.
			Substituting  $\widehat{n-3}= \widehat{n-2}= \varepsilon$ into~\eqref{eq:lem_zn_not_appear}, we get
			$\varphi_m(z_{n-1}) \,\prec \,\varphi_m(P_{n-1}).$
			Part (2) of Lemma~\ref{Lem_start_z_m} implies that $z_{n-1}$ ends with $\underline{n-3}$. 
			Thus, using Part (2) of Lemma~\ref{Lem_varphi_factor},
			$z_{n-1}$ is a factor of $P_{n-1}$, since $m\nmid \,n-2$ implies $\underline{n-2} \neq m$ and thus $\underline{n-3}\neq m-1$, which contradicts the induction hypothesis.
		
			\textbf{Case 2}. Assume that $m \mid n-3$ and $m \nmid n-2$. Replacing $\widehat{n-3}= 0$ and  $\widehat{n-2}= \varepsilon$ into~\eqref{eq:lem_zn_not_appear}, we find  that
				$0^{-1}\,\varphi_m(z_{n-1}) \,\prec \,\varphi_m(P_{n-1})0.$
			By Part (2) of Lemma~\ref{Lem_start_z_m}, $z_{n-1}$ begins and ends with 0 and $P_{n-1}$ begins with $z_0=0$.  Thus $0^{-1}\varphi_m(z_{n-1})$ (resp., $\varphi_m(P_{n-1})0$) begins and ends with 1 (resp., 0).  It follows that $\varphi_m(z_{n-1})$ is a factor of $\varphi_m(P_{n-1})$. As in Case (1), we reach a contradiction since $\underline{n-3} \neq m-1$.	
			
			\textbf{Case 3}. Assume that $m \nmid n-3$ and $m \mid n-2$. Plugging $\widehat{n-3}= \varepsilon$ and  $\widehat{n-2}= 0$ into~\eqref{eq:lem_zn_not_appear}, we have
			$\varphi_m(z_{n-1}) \, 0 \,\prec \,\varphi_m(P_{n-1}).$
			So, we have 	$\,\varphi_m(z_{n-1}) \, 0\, \prec \,\varphi_m(P_{n-1})\,0.$ From Part (1) of Lemma~\ref{Lem_varphi_factor}, 
			$z_{n-1} \,\prec\, P_{n-1}$ contradicting the induction hypothesis.
			\qedhere	
		\end{proof}
		\begin{lemma}
			\label{Prop_closed_z_Tribo}
			For all $n \geq 0$, $z_n$ is closed.
		\end{lemma}
		\begin{proof}
			The result follows from Lemma~\ref{Lemma_w_i_closed} for the case $m=2$. Now suppose that $m \geq 3$.
			Using Table~\ref{table:first-few-words-z-n}, the result is clearly true for the cases $n \in \{0, 1,2\}$. We assume that $n \geq 3$. Let us show, equivalently, that there exists a border $v$ of $z_n$ which $|z_n|_{v}=2$, that is, we prove that $z_n$ is a complete return to $v$. Set $v=(\underline{n-3})^{-1}\, h_{n-3}^R\;\underline{n-2}$. Clearly, $v$ is a border of $z_n$. 
			There are two cases to consider according to the value of $n$.

				\textbf{Case 1}. Suppose that $n \leq m-1$.
				Observe that using Remark~\ref{Rem_z_m},
				$$ z_n = (n-3)^{-1}\,\, h_{n-3}^R \,\, h_{n-2}^R\,n\,	\prod_{i=0}^{n-3} h_i^R \,(n-2).$$		
				In order to show that $z_n$ is a complete return to $v$, it suffices to prove that $y$ is a return word of $v$, where 
				\begin{align}
					y&=(n-3)^{-1}\, h_{n-3}^R \; h_{n-2}^R\;\underline{n}\,\prod_{i=0}^{n-4} h_i^R \,(n-3)\label{formula_y1}.
				\end{align}
				For this purpose, we prove that $fyf^{-1}$ is a return word of $u_{j(v)}$, where $f$  and ${j(v)}$ are those of Lemma~\ref{Lem_Justin}. First, we find the minimal integer $j(v)$ such that the word $v$ is a factor of $u_{j(v)}$.
				We obtain from Lemma~\ref{Lem_t_accurred_in_u_} that $v$ does not occur in $u_{n-1}$.
				From Part (4) of Lemma~\ref{u-h-relations}, we have
				$h_{n-3}^R\; h_{n-2}^R \rhd u_n$.
				Using Lemma~\ref{lem:starthR}, $h_{n-3}^R$ (resp., $h_{n-2}^R$) starts with $n-3$ (resp., $n-2$). Thus,
				$(n-3)^{-1}\, h_{n-3}^R\;(n-2)\,\prec\, u_n$.
				We conclude that $j(v)=n$.	Now,  let
				\begin{align}
					f&=\prod_{i=0}^{n-4} h_i^R\;(n-3).\label{formula_f1}
				\end{align}
				Plugging~\eqref{formula_f1} and~\eqref{formula_y1}  into $fyf^{-1}$ and then using Part (4) of Lemma~\ref{u-h-relations}, we obtain
				\begin{align*}
					fyf^{-1} 
					=\prod_{i=0}^{n-2} h_i^R\;n\;
					=u_{n}\;n.
				\end{align*}
				Using Part (2) of Lemma~\ref{Pn}, we know that
				$
				fyf^{-1}=\mu_{n-1}\,(n).$
				Then, by Lemma~\ref{Lem_Justin}, $y$ is a return word of $v$. So, the desired conclusion is obtained in this case.
				
				\textbf{Case 2}. Suppose that $n\geq m$. The proof is obtained in the same manner as the first case.
				By Remark~\ref{Rem_z_m}, we find that
				$$z_n=(\underline{n-3})^{-1}\, h_{n-3}^R \, h_{n-2}^R\, \prod_{i=(n-m)_*}^{n-3}h_{i}^R\, \underline{n-2}.$$
				Now let
				\begin{align}
					y&=(\underline{n-3})^{-1}\, h_{n-3}^R \; h_{n-2}^R\; \prod_{i=(n-m)_*}^{n-4}h_{i}^R\;\underline{n-3}.\label{formula_y}
				\end{align}
				Similarly to the first case,  the minimal integer ${j(v)}$ such that $v$ is a factor of $u_{j(v)}$ equals $j(v)=n$.
				Substituting~\eqref{formula_f1} and~\eqref{formula_y} into $fyf^{-1}$, we get
				\begin{align*}
					fyf^{-1} 
					&=\prod_{i=0}^{n-4}h_{i}^R\; h_{n-3}^R \; h_{n-2}^R \,\prod_{i=(n-m)_*}^{n-4}h_{i}^R\, (\prod_{i=0}^{n-4}h_{i}^R)^{-1} 
					= \prod_{i=0}^{n-2}h_{i}^R \; (\prod_{i=0}^{n-m-1}h_{i}^R)^{-1}.
				\end{align*}
				By Part (4) of Lemma~\ref{u-h-relations}, we have
				$
				\label{eq:u}
				fyf^{-1}=u_n\, u_{n-m+1}^{-1} \,= \, h_{n-2}\, \cdots \, h_{n-m}.
				$
				Now using  Part (1) of Lemma~\ref{Pn}, we find that
				$
				fyf^{-1}=\mu_{n-1}(\underline{n-m}).
				$
				We conclude from Lemma~\ref{Lem_Justin} that $y$ is a return word of $v$. So, we obtain the desired conclusion in this case.\qedhere
		\end{proof}
		\begin{lemma}
			\label{Lemma_t_factor_u}
			For all $n \geq 4$, the word $(\underline{n-3})^{-1}h_{n-3}^R$ is not a factor of $u_{n-3}$.
		\end{lemma}
		\begin{proof}
			The proof is by induction on $n$. The case $n=4$ can  be easily checked by hand. Assume that $n\geq 5$ and  that the claim holds true for all values less than $n$ and consider the case $n$. Suppose to the contrary that 
			$(\underline{n-3})^{-1}\, h_{n-3}^R \prec \,u_{n-3}.$
			Using Parts (3) and (4) of Lemma~\ref{Pn}, 
			\begin{equation*}
				(\underline{n-3})^{-1}\,0^{-1}\,\varphi_m( h_{n-4}^R)\,0\, \prec \, \varphi_m(u_{n-4})\,0 .
			\end{equation*}
			So, 
			$
			(\underline{n-3})^{-1}\,0^{-1}\,\varphi_m\big( (\underline{n-4})\,(\underline{n-4})^{-1} h_{n-4}^R\big)\,0\, \prec \, \varphi_m(u_{n-4})\,0.
			$
			Using Lemma~\ref{varphi}, we get
			\begin{equation}
				\label{eq:lem:t_factor_u}
				(\widehat{n-3})^{-1}\,\varphi_m\big( \,(\underline{n-4})^{-1} h_{n-4}^R\big)\,0\, \prec \, \varphi_m(u_{n-4})\,0.
			\end{equation}
			Now we divide the proof into two cases according to whether $m$ divides $n-3$ or not.

				\textbf{Case 1}. Suppose that $m  \nmid \,(n-3)$. 
				We obtain that
				$\varphi_m( \,(\underline{n-4})^{-1} h_{n-4}^R)\,0\, \prec \, \varphi_m(u_{n-4})\,0$
				by replacing
				$\widehat{n-3}=\varepsilon$ into \eqref{eq:lem:t_factor_u}.
				Part (1) of Lemma~\ref{Lem_varphi_factor} implies that
				$(\underline{n-4})^{-1}\, h_{n-4}^R \, \prec \, u_{n-4}\,$ which contradicts the induction hypothesis.
				
				\textbf{Case 2}.  Suppose that $m  \mid \,(n-3)$.
				we find that
				$0^{-1}\,\varphi_m\big( \,(\underline{n-4})^{-1} h_{n-4}^R\big)\,0\, \prec \, \varphi_m(u_{n-4})\,0$ by plugging
				$\widehat{n-3}=0$ into~\eqref{eq:lem:t_factor_u}.
				Using Lemma~\ref{lem:starthR}, $0^{-1}\,\varphi_m\big( \,(\underline{n-4})^{-1} h_{n-4}^R\big)$ begins with $1$. Also,  from Lemma~\ref{lem:starthR} and Part (4) of Lemma~\ref{u-h-relations}, $\varphi_m(u_{n-4})$ begins with $0$. Therefore,
				$\varphi_m\big( \,(\underline{n-4})^{-1}h_{n-4}^R\big)\,0\, \prec \, \varphi_m(u_{n-4})\,0.$
				As in Case (1), we reach a contradiction. \qedhere
		\end{proof}
		
		In the following theorem, we obtain the closed $z$-factorization of the $m$-bonacci word based on the sequence of words $z_n$.
		\begin{theorem}
			The closed $z$-factorization of the $m$-bonacci word is
			\begin{equation*}
				cz(h_{\omega}) = (z_0 \, , \, z_1 \, , z_2 \, , \, z_3 \,  ,\, \ldots).
			\end{equation*}
		\end{theorem}
		\begin{proof}
			First note that using Theorem~\ref{ThmFibo}, the case $m=2$ is covered.  Suppose that $m \geq 3$.
				From Lemma~\ref{Lem_factorization_z_m}, $h_{\omega} \,=\, \prod_{n \geq 0} {z}_{n}$. So,  the first few factors of $h_{\omega}$ are $\{z_0, \ldots z_6\}$. It can be easily checked that $z_i$, $1 \leq i \leq 6$, are closed $z$-factors of $h_{\omega}$. Now assume that $n \geq 7$.	
			In order to prove the statement, we need to show the following three claims.
			\begin{description}
				\item[1.]  The word $z_n$ is closed.
				\item[2.]  The word $z_n$ does not appear in ${P}_{n+1}\,\gamma_n^{-1}$ where $\gamma_n$ is the last letter of $z_n$.
				\item[3.]  Every closed prefix of $z_n$ has already appeared in ${P}_n$.	
			\end{description}	
			Claim (1) is true by Lemma~\ref{Prop_closed_z_Tribo} and Claim (2) is true respectively by Lemmas~\ref{Lem_zn_notAppear_zn+1}, ~\ref{Lem_zn_notAppearBefore} and ~\ref{Prop_z_Not_Appear_Before}. To prove Claim (3), we find the largest closed prefix of $z_n$ and prove that this prefix has already appeared in $P_n$. 
			Lemma~\ref{Prop_closed_z_Tribo} implies that a frontier of $z_n$ is $(\underline{n-3})^{-1}\,h_{n-3}^R\,\underline{n-4}$. Now we set $v=(\underline{n-3})^{-1}\, h_{n-3}^R$ and find the closed prefix of $z_n$ with  border $v$.
			We divide the proof into two cases according to the value of $n$.
			
			\textbf{Case 1}. Assume that $n \leq m-1$. In order to prove the statement, we prove that the word $y$ is a return word of $v$, where
				\begin{equation}
					y=(\underline{n-3})^{-1}\,h_{n-3}^R\;\underline{n-2}.\label{formula_y_c1}
				\end{equation}
				To show this, we need to prove that $fyf^{-1}$ is a return word of $u_{j(v)}$ where $f$  and ${j(v)}$ are those of Lemma~\ref{Lem_Justin}.  Using Part (5) of Lemma~\ref{u-h-relations},  
				$
				u_{n-2}=(\underline{n-3})^{-1}\, h_{n-3}^R.
				$
				Therefore, the minimal integer ${j(v)}$ such that  $v$ is a factor of $u_{j(v)}$ equals $j(v)=n-2$ and thus, set $f=\varepsilon$.  Replacing $f=\varepsilon$ and Equation~\eqref{formula_y_c1} into $fyf^{-1}$, we get
				$
				fyf^{-1}
				= (\underline{n-3})^{-1}\, h_{n-3}^R (n-2).
				$
				By Part (5) of Lemma~\ref{u-h-relations} and then Part (2) of Lemma~\ref{Pn}, we find that
				\begin{equation*}
					fyf^{-1} =u_{n-2}(n-2)=\mu_{n-3}(n-2).
				\end{equation*}
				Lemma~\ref{Lem_Justin} implies that  $y$ is a return word of $v$. So,  $(\underline{n-3})^{-1} h_{n-3}^R\,(\underline{n-2})(\underline{n-3})^{-1}\,h_{n-3}^R $ is a closed prefix of $z_n$.		
				Using  Definition~\ref{Def_z_m} and replacing $z_{n-1}$, $z_{n-2}$ and $z_{n-3}$ into  $P_n$,
				 we find that $(\underline{n-3})^{-1} h_{n-3}^R\,(\underline{n-2})(\underline{n-3})^{-1}\,h_{n-3}^R $ is a factor of $P_n$ which  ends the proof in this case.
				 
				\textbf{Case 2}. Suppose that $n \geq m$. This case is similar to the first case.  Set
				\begin{equation}
					y=(\underline{n-3})^{-1}\,h_{n-3}^R\,(h_{n-m-3}^R)^{-1}\;\underline{n-3}.\label{formula_y3}
				\end{equation}
				Lemma~\ref{Lemma_t_factor_u}  tells us that $v$ is not a factor of $u_{n-3}$. On the other hand, using Part (4) of Lemma~\ref{u-h-relations}, we get
				\begin{equation*}
					u_{n-2}=h_0^R\;\cdots\;h_{n-4}^R\, 
					= h_0^R\;\cdots\; h_{n-m-4}^R\, h_{n-3}^R
					= {h_0}^R\;\cdots\; h_{n-m-4}^R\,\underline{n-3}\, (\underline{n-3})^{-1}\, h_{n-3}^R.
				\end{equation*}
				Therefore, the minimal integer ${j(v)}$ such that  $v$ occurs in $u_{j(v)}$ equals $j(v)=n-2$.	 Now set
				\begin{align}
					f&={h_0}^R\,{h_1}^R\;\cdots\; h_{n-m-4}^R\,\underline{n-3}.\label{formula_f3}
				\end{align}
				Substituting  Equations~\eqref{formula_y3} and~\eqref{formula_f3} into $fyf^{-1}$, we find that
				$$
				fyf^{-1}
				= h_0^R\,h_1^R\;\cdots\; h_{n-4}^R\;(h_0^R\,h_1^R\;\cdots\; h_{n-m-3}^R)^{-1}.
				$$
				Using Part (4) of Lemma~\ref{u-h-relations} and then Part (1) of Lemma~\ref{Pn}, we have
				\begin{align*}
					fyf^{-1}
					&= u_{n-2}\,u_{n-m-1}^{-1}=h_{n-4}\,h_{n-5}\, \cdots \, h_{n-m-2}
					=\mu_{n-3}(\underline{ n-m-2}).
				\end{align*}
				Lemma~\ref{Lem_Justin} implies that $ fyf^{-1}$ is a return word of $u_{n-2}$. Thus, $y$ is a return word of $v$, that is, the word $(\underline{n-3})^{-1} h_{n-3}^R\,(h_{n-m-3}^R)^{-1}\,h_{n-3}^R $ is a closed prefix of $z_n$.
				The result is obtained by using  Definition~\ref{Def_z_m} and substituting $z_{n-1}$, $z_{n-2}$ and $z_{n-3}$ into  $P_n$. \qedhere			
		\end{proof}
		
		\section{ Relation between the palindromic and closed\\ $z$-factorizations of the $m$-bonacci words}
		In this section, we link two kinds of factorizations of the $m$-bonacci words, namely the palindromic and closed $z$-factorizations. In \cite{jahannia}, we introduced a variation of the $z$-factorization, the palindromic $z$-factorization, in which each factor is  palindromic. Also, we computed this factorization for the Fibonacci word and more generally for the $m$-bonacci words. The palindromic $z$-factorization of a word $\boldsymbol{w}$ is $pz(\boldsymbol{w})=(p_1\,  , p_2 \,  , \ldots)$ such that $p_i$ is the shortest palindromic prefix of $p_i\,p_{i+1} \cdots$ which occurs exactly once in $p_1\,p_2\,\cdots\, p_i$.
		
		In the following lemma, the length of the $n$-th palindromic $z$-factor of the $m$-bonacci word is expressed by the previous $m$ palindromic $z$-factors.
		\begin{lemma}
		$\mathrm{	\cite[Corollary \,26]{jahannia}}$
			\label{cor:comparison-of-length-mbonacci-2}
			Let $pz(h_{\omega})=(p_1\,  , p_2 \,  , \ldots)$ be the palindromic $z$-factorization of the $m$-bonacci word.
			If $m$ is even, then, for all $n\ge m-1$, we have
			$$
			|p_n| =
			|p_{n-1}| + |p_{n-2}| + \cdots + |p_{n-m}|.
			$$
			If $m$ is odd, then, for all $n\ge m-1$, we have
			$$
			|p_n| = |p_{n-1}| + |p_{n-2}| + \cdots + |p_{n-m}| + (-1)^n.
			$$
		\end{lemma}
		We compare these two types of factorizations in the following corollary.
		\begin{corollary}Let $pz(h_\omega)=(p_1,p_2, \ldots)$ and $cz(h_\omega)=(z_1,z_2,\ldots)$ be respectively the palindromic $z$-factorization and the closed $z$-factorization of $h_{\omega}$.
		If $m=2$, then for all $n \geq m-1$, $|z_n|=|p_n|$.
		If $m\geq 3$, then for every even integer $m$ and for all $n \geq m-1$,  $|z_n|=|p_n|$.
		\end{corollary}
		\begin{proof}
		The case $m=2$ follows from Theorem~\ref{ThmFibo} and the case $m\geq 3$ follows from  Part (1) of Lemma~\ref{Lem_start_z_m} and Lemma~\ref{cor:comparison-of-length-mbonacci-2}.
		\end{proof}
		\section{The oc-sequence of the $m$-bonacci words}
		The notion of the oc-sequence of a word $w$ is introduced in \cite{Luca}. It is a binary sequence whose $n$-th element is 1 if the length-$n$ prefix of $w$ is closed; otherwise, it is 0. In this section, our aim is to show that the sequence of the lengths of the maximum consecutive 1's in $oc(h_{\omega}^{(m)})$ is exactly the $m$-bonacci word.
		
		In the following lemma, we prove that the palindromic prefixes of the $m$-bonacci words are closed.
		\begin{lemma}
			\label{Lem_u_n_Closed}
			For all $n \geq 2$, $u_n$ is closed.
		\end{lemma}
		\begin{proof}
			To prove this statement, it suffices to find a border $v$ of $u_n$ such that $|u_n|_v=2$.
			Using Lemma~\ref{u-h-relations}, we have
			$u_n\,=\,h_{n-2}\, u_{n-1}$. As $u_n$ is a palindrome by its definition, we find that $u_{n-1}$ is a border of $u_n$. We show that $v=u_{n-1}$ .	
			We proceed by induction on $n$. The result is clear for $n\in \{ 2,3\}$ by Table~\ref{table:first-few-words-u-n}. Now suppose that $n\geq 4$ and the result holds true up to $n-1$. We will show that it is still true for $n$. We proceed by contradiction and suppose that $u_{n-1}$ is a proper factor of $u_n$. Therefore, there exist non-empty words $s$ and $t$ such that $u_{n}=s\,u_{n-1}\,t$. 	
			From Part (4) of Lemma~\ref{Pn}, we obtain that
			$
			\varphi_m(u_{n-1})0 =s \,\varphi_m(u_{n-2})0 t.
			$
			From Lemma~\ref{lem:code}, we deduce that $\varphi_m(u_{n-1})$ and $\varphi_m(u_{n-2})$ have a unique factorization using words of the set $C=\{01,02,\ldots,0(m-1),0\}$. Thus, we have $\varphi_m(u_{n-1})0=y_1 \cdots y_k0$ and $\varphi_m(u_{n-2})=x_1\cdots x_{\ell}$, with $x_1,\ldots,x_{\ell},y_1,\ldots , y_k \in C$ and $\ell,k \geq 1$. 
			By uniqueness of the factorization, there exists   $1 \leq i  \leq k$ such that for all $1\leq j \leq l$, we have $x_j= y_{i+j-1}$. Also, $s=y_1\cdots y_{i-1}$ and $0t= y_{i+l}\cdots y_k0$. Thus, there exist words $s',t' \in A_m^{*}$ such that $\varphi_m(s')=s$ and $\varphi_m(t')=t$. Finally, we deduce that
			$$
			\varphi_m(u_{n-1})0=s\varphi_m(u_{n-2})0t=\varphi_m(s')\varphi_m(u_{n-2}) \varphi_m(t')=\varphi_m(s'u_{n-2}t').
			$$
			By injectivity of $\varphi_m$, $u_{n-2}$ is a proper factor of $u_{n-1}$, contradicting the induction hypothesis.
		\end{proof}
		\begin{definition}
			\label{t-n}
			Let $n\ge 2$. We define $t_n$ by
			\begin{equation*}
				t_n =
				\begin{cases}
					(n-1)h_0^R\cdots h_{n-3}^R, &\text{ if   } 2 \leq n \leq m-1;\\
					h_{n-m-1}^R h_{n-m}^R \cdots h_{n-3}^R,& \text{if   } n \geq m. 
				\end{cases}
			\end{equation*}
		\end{definition}
		Then we have, $h_{n-1}^R= t_n h_{n-2}^R$. Hence, $u_n t_n h_{n-2}^R=u_{n+1}$.
		\begin{definition}
			\label{prefix-types}
			Let $w$ be a prefix of the $m$-bonacci word $h_{\omega}$ and let $n(w)$ be the unique positive integer satisfying $|u_{n(w)}|<|w|\leq |u_{n(w)+1}|$.	Then, $w$ is a prefix of {\it type-1} if it satisfies
				\begin{align}
				|u_n(w)|+|t_{n(w)}|<|w|\leq |u_{n(w)+1}|\label{ineq1}
			\end{align}
			and of {\it type-2} if it satisfies 
			\begin{align}
				|u_{n(w)}|<|w|\leq |u_{n(w)}|+|t_{n(w)}|. \label{ineq2}
			\end{align}
		\end{definition}		
		\begin{lemma}
			\label{Two-type-pref}
			Let $w$ be a prefix of the $m$-bonacci word $h_{\omega}$ satisfying $|u_{n(w)}|<|w|\leq |u_{{n(w)}+1}|$ for some positive integer $n(w)$ and let $v=u_{n(w)}^{-1}w$.
			\begin{enumerate}
				\item If $w$ is a prefix of type-1, then $v=t_{n(w)} x$, where $x$ is a non-empty prefix of $h_{{n(w)}-2}^R$.
				\item If $w$ is a prefix of type-2, then $v$ is a non-empty prefix of $t_{n(w)}$ (possibly $v=t_{n(w)}$).
			\end{enumerate}
		\end{lemma}
		\begin{proof}
			The result directly follows from Inequalities~\eqref{ineq1} and~\eqref{ineq2} since $u_{{n(w)}+1} = u_{n(w)} t_{n(w)} h_{{n(w)}-2}^R$ and $w=u_{n(w)}v$.
		\end{proof}

		In the following theorem, we characterize closed (and open) prefixes of the $m$-bonacci word.
		\begin{theorem} We have the following properties.
			\begin{enumerate}
				\item
				The prefixes of type-1 of the $m$-bonacci word are closed.		
				\item
				The prefixes of type-2 of the $m$-bonacci word are open.
			\end{enumerate}
		\end{theorem}
		
		\begin{proof}
			First observe that $0=u_2$ is the type-1 prefix of $h_{\omega}$. Now we examine the prefixes of the $m$-bonacci word whose length is greater than 1. 
			To prove the first item in the statement, consider each prefix $w$ of $h_{\omega}$ with $w=u_{n(w)} v$ and $v=t_{n(w)} x$ such that $\varepsilon \neq x \lhd h_{{n(w)}-2}^R$.
			In order to prove the assertion, we need to find a frontier $u'$ of $w$.
			There are two cases to consider according to the value of ${n(w)}$.

				\textbf{Case 1}. Assume that ${n(w)}\leq m-1$. Therefore, $t_{n(w)}=({n(w)}-1)h_0^{R}\cdots h_{{n(w)}-3}^{R}$. As ${n(w)}-1$ does not occur in $u_{{n(w)}-1}$, the longest border of $w$ is $h_{0}^R h_{1}^R\cdots h_{{n(w)}-3}^Rx= u_{{n(w)}-1}x$.  Using the proof of Lemma~\ref{Lem_u_n_Closed}, $u_{{n(w)}-1}$ has no internal occurrence in $u_{n(w)}$. Hence, $u_{{n(w)}-1}x$ has no internal occurrence in $w$ and we are done.
				
				\textbf{Case 2}. Assume that ${n(w)}\geq m$. Thus, $t_{n(w)}=h_{{n(w)}-m-1}^R h_{{n(w)}-m}^R \cdots h_{{n(w)}-3}^R$. Set $u'=u_{{n(w)}-m}t_{n(w)}x=u_{{n(w)}-1}x$.  It is obvious that $u'$ is a suffix of $w$. Also by Part (3) of Lemma~\ref{u-h-relations}, $u'$ is a prefix of $w$. As in Case (1), $u_{{n(w)}-1}x$ has no internal occurrence in $w$ and we get the result.			
			To prove the second item in the statement, let $w=u_{n(w)} v$, where $v \lhd t_{n(w)}$.
			We consider two cases according to the value of ${n(w)}$.

				\textbf{Case 1}. Suppose that ${n(w)}\leq m$. Then, we have $t_{n(w)}=({n(w)}-1)h_0^{R}\cdots h_{{n(w)}-3}^{R}$. Since ${n(w)}-1$ does not occur in $u_{n(w)}$, the longest border of $w$ is $({n(w)}-1)^{-1}v$ which is also a border of $u_{n(w)}$. So, it appears three times in $w$. In other words, $w$ is open.
				
				\textbf{Case 2}. Suppose that ${n(w)}\geq m$.  We prove that the longest border of $w$ is $u_{{n(w)}-m}\,v$.
				Suppose that the word $xv$ is a border of $w$. As $w=u_{n(w)}v$, $x$ is a border of $u_{n(w)}$. Since $u_{n(w)}$ is palindromic, $x$ is a palindromic prefix of $u_{n(w)}$. Thus it is of the form $x=u_i$ for some $i$.
				As $v$ is a prefix of $t_{n(w)}=h_{{n(w)}-m-1}^R h_{{n(w)}-m}^R \cdots h_{{n(w)}-3}^R$, the longest $u_i$, $i\geq 2$, such that $u_iv$ is a border of $w$ is $u_{{n(w)}-m}$.
				On the other hand, $u_{{n(w)}-m}\,v$ is a border of $u_{n(w)}$.
				Therefore, it occurs three times in $w$ and thus $w$ is open.\qedhere
		\end{proof}
		\begin{corollary}
			The sequence of lengths of the maximum consecutive 1's in the oc-sequence  of the $m$-bonacci word is exactly $|h_i|$, for all $i \geq 0$.
		\end{corollary}
		\begin{proof}
			Using Inequality \eqref{ineq1}, we deduce that the number of consecutive closed prefixes of the $m$-bonacci words $h_{\omega}$ is equal to $|t_n|+|h_{n-2}^R|=|h_{n-1}^{R}|$ and we get the result.
		\end{proof}

		\begin{example}
			Let $t_n$ be the sequence of Tribonacci numbers where $T_{-1}=1$, $T_0=1$, $T_1=2$ and for all $n\geq 2$, $T_n=T_{n-1}+T_{n-2}+T_{n-3}$. Then, we have $oc(h_\omega^{(3)})= 1\,0\,\prod_{i \geq 0} 1^{T_i}\, 0^{T_{i-1}+T_{i}}  $.	
			Table~\ref{table:first-few-words-oc-Tribo} shows the first few values of the $oc$-sequence for the infinite Tribonacci word.
			
			\begin{table}[h]
				\centering	
				\resizebox{\textwidth}{!}
				{
				\begin{tabular}{c | c c c c c c c c c c c c c c c c c c c c c c c c}
						$n$&  1 & 2 & 3 & 4 & 5 & 6 & 7 & 8 & 9 & 10& 11&12&13&14&15&16&17&18&19&20&21&22&23&24\\
						\cline{1-25}
						$h_\omega^{(3)}$&0&1&0&2& 0& 1& 0& 0& 1& 0& 2& 0& 1& 0& 1& 0& 2& 0& 1& 0&0&1&0&2\\
						\cline{1-25}
						$oc(h_\omega^{(3)})$& 1 & 0 & 1 & 0 & 0 & 1 &  1& 0 & 0 & 0 & 1&1&1&1&0&0&0&0&0&0&1&1&1&1\\
					\end{tabular}
				}
				\caption{The first few values of the $oc$-sequence of the Tribonacci word}
				\label{table:first-few-words-oc-Tribo}
			\end{table}
		\end{example}
		\section{Open problems}
		It is interesting to find the closed $z$-factorization of other infinite words such as episturmian words and automatic words. Another interesting problem is to obtain the closed $c$-factorizations of infinite words and find a relation between the closed $z$-factorization and the closed $c$-factorization of infinite words. We  leave it as an open problem to characterize the closed $c$-factorization of the $m$-bonacci word.
		\begin{problem}
			Let $c(h_{\omega})=(c_0, \, c_1,\, c_2,\, \ldots )$ be the closed $c$-factorization of the $m$-bonacci word. For all $m \geq 3$ and $n \geq 2m-1$, we conjecture that $|c_i| =|h^{{m}}_{i-m+1}|$.
		\end{problem}
		\section*{Acknowledgements}
		
		Manon Stipulanti is supported by the FNRS Research Grant 1.B.397.20.

		\vskip 0.4 true cm

		\vskip 0.4 true cm
		

	\end{document}